\documentclass[]{article}

\usepackage{mathtools}
\usepackage{amsthm, amssymb}
\usepackage{bm}
\usepackage{graphicx}
\graphicspath{{plots/}}
\usepackage{xcolor}
\usepackage[margin=1in]{geometry}
\usepackage[shortlabels]{enumitem}
\usepackage[english]{babel}
\usepackage{tabularray}
\UseTblrLibrary{siunitx}
\usepackage{csquotes}
\usepackage{titlecaps}
\usepackage[backend=biber, style=numeric-comp,
sorting = none,
giveninits,
maxbibnames=99]{biblatex}
\Addlcwords{a in of the for an or and by to with}
\DeclareFieldFormat{titlecase}{\titlecap{#1}}
\usepackage{hyperref}

\newcommand*{\R}{\mathbb R}
\newcommand*{\Rp}{{\R_+}}
\newcommand*{\cId}{\mathbb I}
\newcommand*{\J}{\mathcal J}

\renewcommand*{\P}{\bm{P}}
\newcommand*{\Pij}{P_{ij}}
\newcommand*{\Puij}{\Pij^u}
\newcommand*{\Pu}{\P^{\bm{u}}}
\newcommand*{\dt}{\frac{d}{dt}}
\newcommand{\brho}{\bm{\rho}}

\newcommand{\bmm}{\bm{m}}

\newcommand{\bzero}{\bm{0}}
\DeclarePairedDelimiter\ave{\langle}{\rangle}
\newcommand{\bA}{\bm{A}}
\newcommand{\bB}{\bm{B}}
\newcommand{\bD}{\bm{D}}
\newcommand{\bI}{\bm{I}}
\newcommand{\cI}{\mathcal{I}}

\newcommand{\Prob}{\operatorname{Prob}}
\newcommand{\bP}{\bm{P}}
\newcommand{\bR}{\bm{R}}
\newcommand{\cR}{\mathcal{R}}
\newcommand \commentout[1] {}

\definecolor{revcolor}{RGB}{223, 0, 84}

\renewcommand*{\(}{\begin{equation}}
\renewcommand*{\)}{\end{equation}}
\DeclarePairedDelimiter{\ev}
                        {\langle}{\rangle}

\DeclareMathOperator*{\argmax}{arg\,max}
\DeclarePairedDelimiter{\abs}{\lvert}{\rvert}

\let\phi\varphi

\newtheorem{theorem}{Theorem}
\newtheorem{lemma}[theorem]{Lemma}
\newtheorem{prop}[theorem]{Proposition}

\theoremstyle{definition}

\theoremstyle{remark}
\newtheorem{remark}[theorem]{Remark}

\title{Optimal control of a kinetic model describing social interactions on a graph}
\author{Jonathan Franceschi\thanks{Department of Mathematics \lq\lq F. Casorati\rq\rq, University of Ferrara, via Machiavelli, 30, Ferrara, Italy, \texttt{jonathan.franceschi@unife.it}} \and Nadia Loy \thanks{Department of Mathematical Sciences ``G. L. Lagrange'', Politecnico di Torino, Corso Duca degli Abruzzi 24, 10129, Torino, Italy, (\texttt{nadia.loy@polito.it}) (Corresponding Author)}}
\date{}

\begin{document}
\maketitle

\begin{abstract}
In this paper, we investigate the optimal control of a kinetic model describing agents that migrate across a graph and interact within its nodes by exchanging a physical quantity. As a prototype application, we focus on the spread of an infectious disease on a graph, where the exchanged quantity corresponds to the viral load. The control is applied separately to both microscopic mechanisms, the mobility  mechanisms and the interaction dynamics respectively, with the objective of minimizing the average viral load across the system.

By analyzing the macroscopic equations derived from the kinetic model, we demonstrate that the most effective and efficient strategy is to minimize the average viral load weighted by the mass (i.e., agent density) at each node. We explore this approach under two distinct interaction models: in the first, infection (gain) and healing (loss) occur within the same interaction; in the second, infection and healing are modeled as two separate, independent processes.

For the first model, we show that it is possible to halt the progression of the disease, albeit at a very high control cost. In contrast, the second model allows for the complete eradication of the disease with a more moderate control effort. Numerical simulations illustrate the role of each type of intervention and the interplay between mobility and interaction control strategies in both models.
\end{abstract}

\section{Introduction}

Interactions happening within---or among---social groups often present traits of heterogeneity: whenever either the groups or the interactions themselves are not homogeneous, labeling the social groups' elements is a useful abstraction, along with building a relational model of members' interactions based on those labels. In this context, network theory is a natural framework~\cite{newman2018networks}, allowing to study different phenomena from disease spreading~\cite{kiss2017mathematics, pastor15} to dissemination of information~\cite{moreno2004dynamics, bovet2019influence} to movement (local, e.g., traffic~\cite{papageorgiou1990dynamic, barthelemy2006optimal} and global, e.g., migration~\cite{taylor2010population, haug2008migration}). In particular, developing models to analyze and to control dynamics taking place on large networks has become one of the most relevant aspects in contemporary research in applied mathematics, especially related to epidemiology.

Kinetic theory~\cite{cercignani1988BOOK}, emerging from the field of statistical mechanics of particle physics, has established itself in the last decades as one of the most powerful frameworks to model the behavior of systems of a large number of interacting agents. Relying on the assumption of indistinguishability of the agents who then behave following a common interaction rule, kinetic theory allows one to derive sound population-level macroscopic models inheriting the characteristics of the implemented individual microscopic dynamics~\cite{pareschi2013BOOK}. For the latter reasons, kinetic models have also been proved to be effective in the study of dynamics on underlying networks, either for diseases spreading~\cite{BertagliaPareschi,loytosin2021}, opinion formation~\cite{albi2014bis,albi2024data,during24} or more general types of interactions~\cite{nurisso24, burger21}.

Regarding the specific application, in the field of epidemiology, viral-load-based frameworks are classically employed, especially in studies on chronic infections like HIV~\cite{kelley07, killingo17}; they are also studied for acute settings in view of the usefulness of quantitative data over qualitative clinical descriptions for policy making or therapy design~\cite{hay21,petrie13}. For this reason, and also for its natural description within kinetic theory, substantial mathematical research on viral-load-based models has been conducted recently, see, e.g.,~\cite{loytosin2020, loytosin2021} for a more general introduction,~\cite{RDMNLAT,RDMNLAT2} for compartmental-based models and~\cite{bondesan24} for a data-driven approach. 

Concerning optimal control theory of diseases on networks, instead, we refer the reader to the recent works on homogeneous populations~\cite{espinoza20, khouzani11, li19}, on heterogeneous populations~\cite{roy11}, to works on information dissemination like, e.g.,~\cite{el20, liu20, kandhway17}, and to works with ODEs settings on networks also considering vaccination~\cite{KABIR2019,kabir2020} and optimal control for compartmental models~\cite{zino1,zino2}. Within the framework of kinetic theory, we mention the following works concerning opinion formation~\cite{albi2014bis} and optimal control on related topics, see, e.g.,~\cite{albi2021, franceschi23, Dimarco2022} and references therein.

In view of the importance of this topic for the aforementioned reasons, in this paper we follow up the model presented in~\cite{loytosin2021}, where the authors propose a Boltzmann-type kinetic model for the spreading of diseases on a graph, based on the exchange of viral load. The authors in~\cite{loytosin2021} formulate  stochastic processes allowing to describe accurately complex microscopic dynamics, such as, e.g., interacting commuters on a graph. From the stochastic processes, kinetic and then population-level macroscopic models are derived. Such models, that are not postulated at the aggregate level, are sound as they include microscopic details of the dynamics. The asymptotic trend of some average population-level quantities is investigated, and the behavior turns out to depend meaningfully on the parameters of the microscopic dynamics. Therefore, our contribution to the model in~\cite{loytosin2021} is to devise optimal control strategies to mitigate the infection on the network. To this aim, we rely on the method used in~\cite{albi2021}.

The manuscript outline is as follows. In Section~\ref{sec:basic} we present the modeling framework and a summary of the results presented in~\cite{loytosin2021} for the sake of the reader. We consider a strongly connected graph: each agent belongs to exactly one of its vertices (or nodes) and can move from node to node with a certain probability, prescribed by a transition matrix. Agents can only interact---pairwise---with their peers belonging to their same node; each binary interaction consists in an exchange of viral load, representing the contagion-healing dynamics. In Section~\ref{sec:control}, we start by considering a control policy intervening on agents' mobility, i.e., we construct a suitably controlled transition matrix. Then we proceed to devise in-node control strategies to affect agents' interactions within the same vertex. Our results show that the prototypical model of~\cite{loytosin2021}, while useful to gain insight by binary exchange contagion dynamics, might be too simplistic to depict a fully realistic infection phenomenon in which individuals can actually recover, and to give rise to satisfactory controlling results in a challenging spreading scenario. For this reason, in Section~\ref{sec:healing} we modify the interaction dynamics to incorporate an infection process, again expressed as a suitable kind of binary interaction, and a healing process, described as an autonomous process. In this new setting, in Section~\ref{sec:healing.control}, we devise a different in-node control strategy, of which we can prove it achieves the complete eradication of the disease from the network, under appropriate assumptions. We present in Section~\ref{sec:num} several numerical experiments to support our theoretical findings and we also showcase an application of the new infection-healing dynamics on real-world data. Finally, we conclude the manuscript by commenting on the presented results and we outline some possible future research directions. For reader's convenience, we report in Table~\ref{tab:symbols} a summary of the most relevant symbols we use.

\begin{table}[htbp]
\centering
\begin{tblr}{columns = {l}, column{1} = {mode=dmath},
hline{1, Z} = {solid, 0.4ex}, vline{1, 3} = {solid, 0.4ex}}
\text{\bfseries Symbol} & \bfseries Definition\\[1ex]
\cI & Set of nodes\\
\P = [\Pij]_{i,j\in\cI} & Transition matrix\\
v & Viral load\\
\nu_1^i,\ \nu_2^i & Node-dependent exchange parameters\\
\chi & Frequency of migration\\
\mu & Frequency of interaction in infection-only model\\
\rho_i & Number of agents/mass of the distribution in node $i$\\
m_i & Average viral load in node $i$\\
u_i^\chi & Control on the migration\\
\Pu = [\Puij]_{i,j\in\cI} & Controlled transition matrix\\
u_i^\mu & In-node interaction control\\
\rho_i^u & Controlled number of agents in node $i$\\
m_i^u & Controlled average viral load in node $i$\\
\delta & Minimum entry of $\Pij$\\
\sigma & Healing frequency\\
\rho_i^c & $\dfrac{\gamma\nu_1^i}{\sigma\nu_2^i}$\\
u_i^\sigma & In-node healing control\\
\mathcal R_0 & Basic reproduction number
\end{tblr}
\caption{Table of symbols.}
\label{tab:symbols}
\end{table}

\section{Mathematical modeling of interactions on a graph}\label{sec:basic}

In this section, we revise the basic kinetic-like model (without control) for social interactions on a graph. We summarize the results presented in~\cite{loytosin2021} concerning the emerging collective behavior of the system; the latter can be studied by analyzing the macroscopic equations for average quantities that are derived from the underlying kinetic description.

\subsection{Kinetic description}

Let us consider a large system of interacting individuals migrating on a network which is modeled by a graph with a finite number of vertices and edges. In particular, we introduce a weighted graph $G=(\cI,E,A)$, where $\cI$ is the set of the \emph{nodes} that is a finite ordered index set with $\abs{\cI}=n\in\mathbb{N}$, for instance $\cI=\{1,\,\dots,\,n\}\subset\mathbb{N}$. The set $E$ is the set of the edges, that is a subset of $\cI\times\cI$, while $A\coloneqq\left(A_{ij}\right)_{i,j\in \cI}$ is the matrix of the weights that can be assigned to each edge ($A_{ij}=0$ corresponds to no edge connecting nodes $i$ and $j$). Then $A$ allows to introduce the \emph{transition matrix} $\P \coloneqq [\Pij]_{i,j\in\cI} \in \R^{n \times n}$ that is defined as 
\(
\Pij \coloneqq \dfrac{A_{ij}}{\sum_{i\in \cI} A_{ij}} \in [0,1], \qquad i,j \in \cI,
\)
and thus satisfies
\begin{equation}\label{eq:P}
\sum_{i\in \cI} \Pij = 1, \qquad \forall j \in \cI.
\end{equation}

The transition matrix $\P$ is, by definition~\eqref{eq:P}, left stochastic.  
Moreover, we shall consider \emph{strongly connected graphs}, which means that for any two nodes, there exists at least one directed path that connects them. Notice that a graph is strongly connected if and only if $\P$ is \emph{irreducible}~\cite{minc1988nonnegative}.

The agents are assumed to be characterized by a physical quantity $v\in \R_+$ and to be located on the nodes of the graph $G$. As a consequence, the agents are characterized by a microscopic state $(x,v)\in \cI\times \R_+$. The label $x\in\cI$, denoting the vertex on which the agent is located, may change as a consequence of possible migrations among the nodes of the graph. Specifically, we define the migration dynamics on the graph as a Markovian process by means of the transition matrix, by defining the probability of migrating from node $j$ to node $i$ as
\begin{equation}\label{def:prob}
\Prob(j \to i)\coloneqq P_{ij}.    
\end{equation}

Conversely, the physical quantity $v\in\R_+$ can be exchanged as a consequence of binary interactions. In the present case $v$ represents the viral load. In particular, it is assumed that individuals may interact only when on the same node $x=i$. In order to have easily manageable interaction rules and having the aim of studying the emerging properties of the system, we consider linear interaction rules defined as
\(\label{eq:update-v}
v' = (1-\nu_1^i+\eta)v + \nu_2^i v_*, \qquad v_*' = (1-\nu_1^i+\eta_\ast)v_* + \nu_2^i v,
\)
where $v$, $v_*$ represent the pre-interaction states of two interacting individuals and $v'$, $v_*'$ their post-interaction states. In~\eqref{eq:update-v},  $\nu_1^i,\nu_2^i\in [0,1]$ are the exchange parameters describing the deterministic part of the interaction that may depend on the node $i$, while $\eta,\eta_\ast$ are white noises taking into account the possible stochastic fluctuations, i.e.
\[
\ave{\eta}=\ave{\eta_\ast}=0, \qquad \ave{\eta^2}=\ave{\eta_*^2}=1,
\]
where, here ad henceforth, $\ave{\cdot}$ denotes expectation.
Moreover, the two microscopic processes leading to the migration across the nodes of the graph and to the binary exchange process are assumed to be stochastically \emph{independent} and happening with frequencies $\chi$ and $\mu$ respectively.

 Discrete-time stochastic processes implementing the prescribed microscopic dynamics~\eqref{def:prob}-\eqref{eq:update-v} can be formulated and allow to derive the  evolution equation of each $f_i=f_i(t,v):\R_+\times\R_+\to\R_+$, which is the statistical  distribution of $v$ in each node $i$. Its collision-like kinetic evolution equation in weak form is 
\(
\begin{aligned}\label{eq:boltz.fi}
\frac{d}{dt}\int_{\Rp} \phi(v) f_i(t,v) \, dv
&= \chi \int_{\Rp} \phi(v) \biggl(\sum_{j \in \cI}\Pij f_j(v,t) -f_i(t,v)\biggr)\, dv\\
&+ \mu \int_{\Rp} \phi(v)Q(f_i,f_i)(v,t)\, dv, \qquad i\in \cI,
\end{aligned}
\)
where $\phi\colon \Rp \to \R$ is a test function and $Q(f_i,f_i)$ is the so-called \emph{collisional operator}, defined in weak form by
\(
\int_{\Rp} \phi(v) Q(f_i,f_i)(v,t)\, dv = \int_{\Rp}\!\int_{\Rp} \ev{\phi(v') - \phi(v)} f_i(t,v)f_i(t,v_*)\, dv_* \,dv. 
\)

\subsection{Aggregate description}\label{sec2:aggr}
In the rest of the manuscript, we will be interested in some aggregate average quantities, namely
\begin{equation}\label{def:macro}
\rho_i(t) \coloneqq \int_{\Rp} f_i(t,v)\, dv, \qquad \rho_i(t)  m_i(t) \coloneqq \int_{\Rp} vf_i(t,v)\, dv,
\end{equation}
which are the mass density and first moment of the individuals on the $i$-th node of the graph, respectively. 

Setting $\varphi=1$ in~\eqref{eq:boltz.fi} allows one to obtain the evolution equation of the $i$th mass, that is
\begin{equation}
	\frac{d\rho_i}{dt}=\chi\left(\sum_{j\in\cI}P_{ij}\rho_j-\rho_i\right), \qquad i\in\cI,
	\label{eq:rho_i}
\end{equation}
which in vector notation reads
\begin{equation}
	\frac{d\brho}{dt}=\chi(\bP-\cId)\brho,
	\label{eq:rho_vect}
\end{equation}
with $\brho\coloneqq(\rho_i)_{i\in\cI}$.
First of all we remark that
$$ \frac{d}{dt}\|\brho(t)\|_1\coloneqq\frac{d}{dt}\sum_{i\in\cI}\rho_i(t)=0, $$
which means that the $\ell^1$-norm of $\brho$, that is the quantity $$\|\brho(t)\|_1\coloneqq\displaystyle\sum_{i\in\cI}\rho_i(t),$$ is conserved in time, which, physically, means conservation of mass across the graph during the whole dynamics.
Specifically, we have that, given $\boldsymbol{\rho^0}=[\rho_1^0,...,\rho^0_n]$ the initial condition, 
\[
\|\brho\|_1= \|\brho^0\|_1.
\]

We choose, without loss of generality, $\|\brho\|_1 \equiv 1$.  
From~\eqref{eq:rho_vect} we can investigate the stationary mass distribution $\brho^\infty\in\R^n_+$ emerging for large times, namely the vector satisfying the equation
\begin{equation}
	(\bP-\cId)\brho^\infty=\bzero.
	\label{eq:rho.inf}
\end{equation}

The study of the eigenvalues-eigenvector properties of the transition matrix $\P$ and application of the Perron-Frobenius theory, allows one to prove that
\begin{prop}{~\cite{loytosin2021}} \label{prop:rho^inf} 
Let $\P$ be irreducible. Then there exists a unique physically admissible solution $\brho^\infty\in\R^n_+$ to~\eqref{eq:rho.inf}, which is a stable and attractive asymptotic density distribution for~\eqref{eq:rho_vect}.
\end{prop}
\begin{remark}
    The uniqueness is fixed by conservation of mass across the graph. 
\end{remark}

Conversely, setting $\varphi(v)=v$ in~\eqref{eq:boltz.fi} allows one to investigate the evolution of the first moment $\rho_i m_i$ 
which turns out to satisfy the equation
\begin{equation}
	\frac{d}{dt}(\rho_im_i)=\chi\left(\sum_{j\in\cI}P_{ij}\rho_jm_j-\rho_im_i\right)+\mu(\nu_2^i-\nu_1^i)\rho_i^2m_i.
	\label{eq:rho_im_i}
\end{equation}
Moreover, we can introduce the average in each node
$$ 
m_i(t)\coloneqq\frac{1}{\rho_i(t)}\int_{\R_+}vf_i(t,v)\,dv, 
$$
so that $\rho_i m_i$ may be seen as the average weighted by the $i$th mass $\rho_i$. Equations~\eqref{eq:rho_im_i} and~\eqref{eq:rho_i} allow one to obtain the evolution equation of $m_i$ that is 
\begin{equation}\label{eq:mi}
	\dt m_i=\frac{\chi}{\rho_i}\sum_{j\in\cI}P_{ij}\rho_j(m_j-m_i)+\mu(\nu^i_2-\nu^i_1)\rho_im_i.
\end{equation}

As expected, the average varies because of the sum of the two independent contributions, which are related to the migration on the graph and to the evolution of the average physical quantity, respectively, as
\begin{equation}
\dt m_i=\left[\dt m_i\right]_\chi +\left[\dt m_i\right]_\mu,
\end{equation}
where
\begin{equation}\label{eq:ave_chi}
\left[\dt m_i\right]_\chi=\frac{\chi}{\rho_i}\sum_{j\in\cI}P_{ij}\rho_j(m_j-m_i),
\end{equation}
and 
\begin{equation}\label{eq:ave_mu}
\left[\dt m_i\right]_\mu=\mu(\nu^i_2-\nu^i_1)\rho_i m_i.
\end{equation}

\begin{figure}[htbp]
    \setbox0=\hbox{\includegraphics[width=0.45\linewidth]{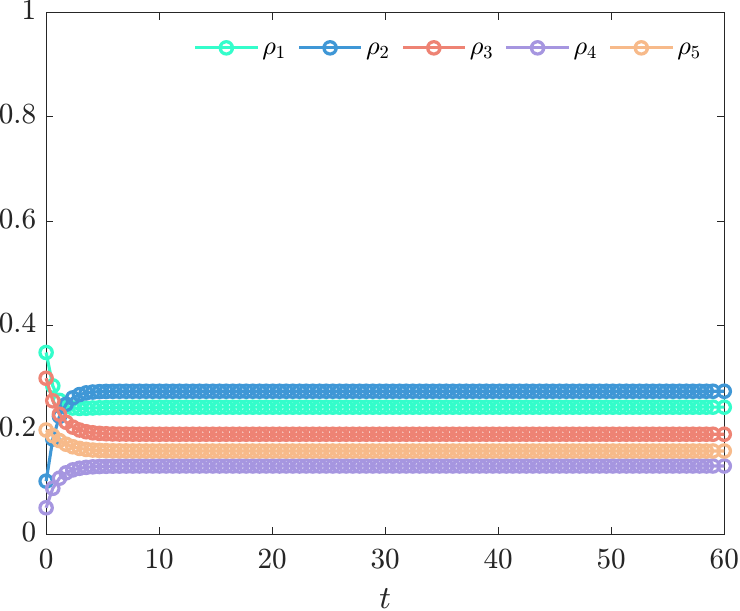}}
    \hbox to \textwidth{%
    \includegraphics[width=0.45\linewidth]{plots/rho-nu1=nu2.pdf}\hfil
    \includegraphics[width=0.45\linewidth, height=\ht0]{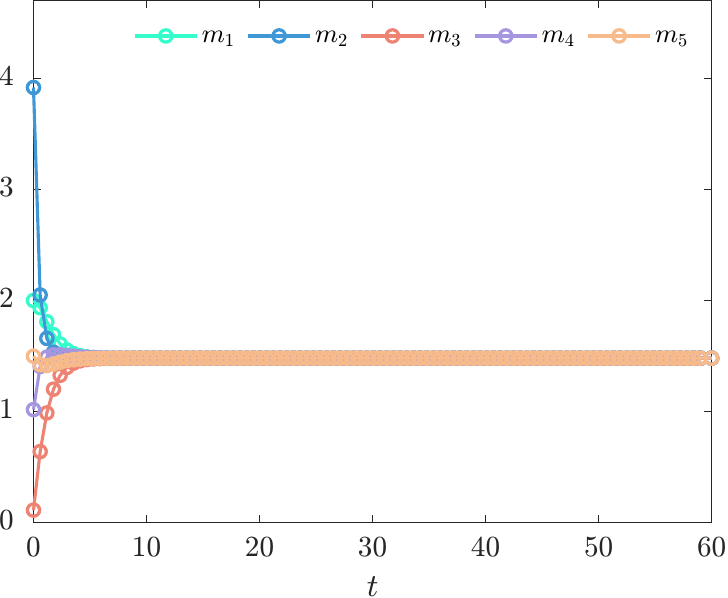}}
    \caption{Left to right: example of evolution in time of number of agents and average viral load in the case $\nu_1^i = \nu_2^i = 1/2$ for all $i \in \cI$. The numerical test shows accordance with Proposition~\ref{prop2}. Refer to Section~\ref{sec:num1} for more details about the simulation.}
    \label{fig:perron-frobenius}
\end{figure}

Moreover, as a consequence of~\eqref{eq:P}, we have that 
\begin{equation}\label{eq:ave_tot}
   \dt \sum_{i\in \cI} \rho_i m_i= \mu \sum_{i\in \cI} \rho_i^2 m_i (\nu_2^i-\nu_1^i),
\end{equation}
which means that the variation of the total weighted average on the graph is only due to the interactions inside the nodes.
As a consequence of the latter, it is possible to prove that 
\begin{prop}\label{prop2}
Assume the graph is strongly connected and that $\nu_1^i=\nu_1, \nu_2^i=\nu_2$. When $t\rightarrow \infty$ the following cases hold true.
\begin{enumerate}
    \item If $\nu_1>\nu_2$ then the solution $m_i$ of~\eqref{eq:mi} satisfies $m_i\rightarrow 0$ for all $i \in \cI$ and $\bmm^{\infty}=\boldsymbol{0}$ is a stable and attractive solution for~\eqref{eq:mi}.
    \item If $\nu_1<\nu_2$ then the solution $m_i$ of~\eqref{eq:mi} satisfies $m_i \rightarrow \infty$ for some $i \in \cI$.
    \item If $\nu_1=\nu_2$ then $\displaystyle\sum_{i\in\cI} \rho_i m_i$, being the evolution of $\rho_i m_i$ described by~\eqref{eq:rho_im_i}, is constant in time. Moreover, there exists a unique stable and attractive equilibrium configuration for $\bmm^\infty$ and in particular $m_i^\infty =  \displaystyle \sum_{j\in\cI} \rho_j m_j$ for all $i \in \cI$. 
\end{enumerate}
\end{prop}
\begin{proof}
The first two points have been proved in~\cite{loytosin2021}. We then prove the third point, as it is slightly more general than in~\cite{loytosin2021}. Let $\nu_1 = \nu_2$. If either $\brho(0) = \bzero$ or $\bmm(0) = \bzero$, equation~\eqref{eq:ave_tot} gives the claim immediately. Otherwise, let us consider 
the quantities $\beta_i$, defined as
\(
\beta_i(t) \coloneqq \frac{\rho_i m_i(t)}{\displaystyle\sum_{j \in \cI} \rho_j m_j(t)}.
\label{eq:betadef}
\)
Since $\nu_1 = \nu_2$, equation~\eqref{eq:ave_tot} implies the conservation of the total first moment $\displaystyle\sum_{i \in \cI}\rho_i m_i$, so that the denominator in equation~\eqref{eq:betadef} is a positive constant.
We can conveniently rewrite equation~\eqref{eq:rho_im_i} as
\(
\frac{d\beta_i}{dt}=\chi\left(\sum_{j\in\cI}P_{ij}\beta_j-\beta_i\right), \qquad i\in\cI,
\label{eq:beta_i}
\)
and~\eqref{eq:beta_i} is well defined. In particular, we also have that $\beta_i \in [0, 1]$ for all $i \in \cI$ and $\sum_{i \in \cI} \beta_i = 1$. Therefore, due to the irreducibility of the matrix $\bP$, we can apply Perron-Frobenius theorem on the vector form of equation~\eqref{eq:beta_i}
\[
\frac{d \bm\beta}{dt} = \chi (\bP - \cId)\bm\beta
\]
to obtain the existence of a unique, stable and attractive positive equilibrium point $\bm\beta^\infty \in \R_+^n$. Moreover, the Perron-Frobenius theorem actually tells us that $\bm\beta^\infty = \brho^\infty$, since they must be scalar multiples and they share the same norm. This implies that
\[
\rho_i^{\infty} m_i^\infty = \rho_i^\infty \sum_{j \in \cI} \rho_j m_j, 
\]
and we conclude by setting the equilibrium point for the node average as
\begin{equation}\label{eq:m_iinfty}
m_i^\infty = 
\begin{cases}
0 & \text{if } \rho_i^\infty = 0,\\
\displaystyle\sum_{j \in \cI} \rho_j m_j & \text{otherwise}.
\end{cases}
\end{equation}
\end{proof}
\begin{remark}
    Proposition~\ref{prop2} allows to detect the asymptotic value of the average viral load in the case $\nu_2=\nu_1$ that is the same in each node as proved by~\eqref{eq:m_iinfty}, and shown numerically in Figure~\ref{fig:perron-frobenius}.
\end{remark}

In the case of node dependent exchange parameters, the authors in~\cite{loytosin2021} prove that
\begin{prop}{\cite{loytosin2021}}\label{prop4} Let the graph be strongly connected and let us assume that $\mu=1/\varepsilon$, $\chi=1$, where $\varepsilon \ll 1$ in~\eqref{eq:boltz.fi}. The following statements hold.
\begin{itemize}
    \item[1] Assume $\nu_1^i>\nu_2^i$ for all $i \in \cI-\lbrace i_\ast\rbrace$ and $\nu_1^{i_*}=\nu_2^{i_*}$. Then $m_i \rightarrow 0$ for all $i \in \cI$, in particular also for $i=i_\ast$.
    \item[2] Assume $\nu_1^i=\nu_2^i$ for all $i \in \cI-\lbrace i_\ast\rbrace$ and $\nu_1^{i_*}>\nu_2^{i_*}$. Let moreover $\cI_\ast=\lbrace j \in \cI: P_{i_\ast,j}>0\rbrace$. Then $m_i \rightarrow 0$ for all $i \in \cI_*$.
    \item[3] Assume $\nu_1^{i_*}<\nu_2^{i_*}$ for some $i_\ast \in \cI$. Then $m_i \rightarrow \infty$ when $t \rightarrow \infty$ for all $i \in \cI$ s.t. $\nu_1^i \le \nu_2^i$, while $m_i \rightarrow 0$ when $t\rightarrow \infty$ for all $i \in \cI$ s.t. $\nu_1^i>\nu_2^i$.
\end{itemize}
\end{prop}
The Proof of the latter (see~\cite{loytosin2021}) relies on the use of classical arguments of kinetic theory, such as the \emph{hydrodynamic regime} and collision invariants. Specifically, the regime $\mu=\dfrac{1}{\varepsilon}$, $\chi=1$ corresponds to the one in which local interactions within the vertices of the graph are much more frequent than jumps from node to node, that are assimilated to the free particle transport in classical kinetic theory, and that happens on a slower time scale. 

\section{The control problem}\label{sec:control}

In this section, we implement a control on the collision-like kinetic equations presented in the previous section.
As a prototype application, we consider the spread of an infectious disease on a graph.
The control has the aim of mitigating the diffusion of the disease on the graph. As the state of the individual with respect to the disease is characterized by the viral-load $v\in \R_+$,  we shall implement such a control by controlling some quantity related to the average viral-load.
As we will show, the most convenient strategy is to control the evolution of the weighted average in each node. In order to do this,  we bear in mind that such an evolution depends on both (independent) processes: the migration across the nodes, and, primarily, the binary interactions leading to the exchange of the physical quantity $v$. As a consequence, we consider two kinds of control: $u_i^\chi$, a multiplicative control term on $\Pij$---with the idea from a modeling point of view of controlling the mobility from node to node---and a multiplicative control $u_i^\mu$ on the binary interactions to mimic containment measures within nodes.

\subsection{Control on the transition matrix}
Let us define the controlled transition matrix $\Pu$ as 
\(
\label{eq:defpuij}
\Puij =
\begin{cases}
\Pij\cdot (1-u_i^\chi) & i \ne j,\\[0.2cm]
1 - \displaystyle\sum_{\substack{k\in \cI\\k\ne i}} P_{ki}^u & \text{otherwise},
\end{cases}
\)
where $u_i^\chi$ is the control on the migration dynamics implemented in order to mitigate the incoming migration in each node $i\in\cI$. We shall consider $u_i^\chi \in [0,1]$, so that $(1-u_i^\chi) \in [0,1
]$ with the idea of expressing migration limitations as a percentage of the normal regime. The optimality conditions defining $u_i^\chi$ will be discussed later, and will be based on the infection condition in node $i$, as it is reasonable to hamper individuals from migrating to node $i$ if the infection there is too high.
The modeling choice $u_i^\chi \in [0,1]$ allows one to prove that the entries of the controlled matrix are non-negative and smaller than one, so that they can be probabilities.

In fact, the following Proposition holds.
\begin{lemma}
If  $u_i^\chi \in [0,1]$ for all $i\in\cI$, then $ \Puij \in [0,1]$ for all $i,j \in \cI$.
\end{lemma}
\begin{proof}
 As $P_{ij} \in [0,1]$ and $u_i^\chi \in [0,1]$ it is clear that $\Puij \in [0,1]$ for all $j \neq i$. Then, we need to check that $P_{ii}^u \in [0,1]$ for all $i \in \cI$. Clearly, $1 - \sum_{k \ne i} P_{ki}^u  \le 1$, while we have that
\[
1 - \displaystyle\sum_{\substack{k\in \cI\\k\ne i}} P_{ki}^u=1 - \displaystyle\sum_{\substack{k\in \cI\\k\ne i}} P_{ki} (1-u_k^\chi) = P_{ii} + \displaystyle\sum_{\substack{k\in \cI\\k\ne i}} P_{ki}\, u_k^\chi \ge 0 ,
\]
where in the last equality we have used $1-\displaystyle\sum_{\substack{k\in \cI\\k\ne i}} P_{ki}=P_{ii}$ that follows from~\eqref{eq:P}.
\end{proof}
Then, as $ \Puij \in [0,1]$ for all $i,j \in \cI$, we can now remark that $\Pu$ is a left stochastic matrix by construction, thanks to the definition of $P_{ii}^u$ in~\eqref{eq:defpuij}, i.e., it satisfies
\begin{equation}\label{Pu_left}
    \displaystyle\sum_{\substack{i\in \cI}} P^u_{ij}=1, \qquad \forall j \in \cI.
\end{equation}
The diagonal entries $P_{ii}^u$ can be explicitly determined from~\eqref{eq:defpuij} using~\eqref{Pu_left} 
\begin{equation}\label{eq:Piiu}
    P_{ii}^u=P_{ii}+\sum_{\substack{k\in \cI\\k\ne i}} u_k P_{ki}.
\end{equation}
The latter means that, while $(1-u_k^\chi)$ diminishes the mobility from node $i$ to $k$, as soon as at least one $k \neq i$ satisfies $u_k^\chi >0$, we have that $P_{ii}^u>P_{ii}$, i.e., the probability of staying in the same node is higher. Conversely, the probability of reaching node $i$ in the controlled and uncontrolled scenarios are related as follows
\[
\displaystyle\sum_{\substack{k\in \cI\\k\ne i}} P^u_{ik}=(1-u_i^\chi) \displaystyle\sum_{\substack{k\in \cI\\k\ne i}} P_{ik}.
\]
The latter means that $(1-u_i^\chi)$ is the total mobility reduction rate to node $i$. Moreover, we have that the control may preserve the irreducibility of the transition matrix under a suitable condition.
\begin{prop}\label{prop:irreducible}
If we set $u_i^\chi \in [0,1)$ for all $i\in\cI$ and consider an irreducible transition matrix $\P$, then $ \Pu$ defined by~\eqref{eq:defpuij} is irreducible.
\end{prop}
\begin{proof}
It is straightforward to see that if $\P$ is associated to a strongly connected graph, then the graph associated to $\Pu$ is still strongly connected and, thus $\Pu$ is irreducible.
\end{proof}
\begin{remark}
Of course, the extremal choice $u_i^\chi = 1$ (that corresponds to a hundred percent restrictions) would have the effect of disconnecting the network, thus ceding the irreducibility. 
\end{remark}

\subsection{Control in each node}
Let us now introduce the control on the kinetic equation~\eqref{eq:boltz.fi}.
We consider two independent control mechanisms on the two independent microscopic processes: the first one on the migration dynamics and the second one on the binary interactions.
We then define the controlled problem as
\begin{equation}\label{eq:boltz.fi.c}
\begin{aligned}
\dt \int_\Rp \phi(v) f_i(t,v)\, dv &= \chi \int_\Rp \phi(v) \left[ \displaystyle\sum_{j\in\cI}\Puij f_j(v,t) - f_i(t,v) \right] \, dv\\
                                          &\hphantom{{}=}+\mu_i\, \ev*{\int_\Rp\!\int_\Rp (\phi(v') - \phi(v)) f_i (v,t) f_i(t,v_*) \,dv_* \, dv} 
\end{aligned}
\end{equation}
where the control matrix $\Puij$ is defined by equation~\eqref{eq:defpuij}, while the control $\mu_i$ on the binary interactions within node $i$ is defined by
\begin{equation}\label{def:mui}
    \mu_i\coloneqq\mu (1-u_i^\mu).
\end{equation}
The latter has the effect of reducing the interaction rate inside each node. The optimality conditions defining $u_i^\mu$ will be discussed later.

Next, we analyze the controlled macroscopic equations for the aggregate quantities. We shall denote by the apex $^u$ the macroscopic quantities related to the controlled problem~\eqref{eq:boltz.fi.c}. Setting $\varphi=1$ in~\eqref{eq:boltz.fi.c}, we obtain
\[
\begin{aligned}\label{eq:rho_u}
\dt \rho_i^u &= \chi \biggl[\sum_{j\in \cI} \Puij \rho_j^u - \rho_i^u \biggr]
                     = \chi \biggl[\sum_{\substack{j\in \cI\\j\ne i}} (\Puij \rho_j^u - P_{ji}^u\rho_i^u )\biggr],\\
                    &= \chi \biggl[(1-u_i^\chi)\sum_{\substack{j\in \cI\\j\ne i}} \Pij \rho_j - \sum_{\substack{j\in \cI\\j\ne i}}P_{ji}(1-u_j^\chi)\rho_i\biggr],
\end{aligned}
\]
that highlights the fact that $(1-u_i^\chi)$ in the gain term regulates the incoming migration, while the terms $(1-u_j^\chi), \, j\ne i$ in the loss term diminish the outgoing flux. 

Setting $\varphi=v$ in~\eqref{eq:boltz.fi.c} we obtain
\(\label{eq:momentum_u}
\dt \rho_i^u(t)m_i^u(t) = \chi \biggl[\sum_{j\in \cI} \Puij \rho_j^u(t)m_j^u(t) - \rho_i^u(t)m_i^u(t) \biggr]
                      + \mu (1-u_i^\mu) (\nu^i_2 - \nu^i_1)\rho_i^{u^2}(t)m_i^u(t),
\)
and, also
\(\label{eq:mean_u}
\begin{aligned}
\dt m_i^u(t)  = \chi \biggl[\sum_{\substack{j\in\cI\\j\neq i}} \Puij \frac{\rho_j^u(t)}{\rho_i^u(t)}\bigl(m_j^u(t) - m_i^u(t)\bigr) \biggr]
                      + \mu (1-u_i^\mu) (\nu^i_2 - \nu^i_1)\rho_i^u(t)m_i^u(t),
\end{aligned}
\)
i.e.
\[
\dt m_i^u=(1-u_i^\chi) \left[\dt m_i\right]_\chi + (1-u_i^\mu) \left[\dt m_i\right]_\mu,
\]
where in the right hand side the quantities are defined in~\eqref{eq:ave_chi}-\eqref{eq:ave_mu}.

\begin{remark} We remark that, for both the evolution of the masses and of the averages, the controls $1-u_i^\chi, 1-u_i^\mu \in [0,1]$ reduce the time variation rate, without inverting the natural trend. Only if we allow $u_i^\mu=u_i^\chi=1$, the control has the effect of stopping the time evolution of both masses and averages. This is why we also consider a control on the migration mechanism.
\end{remark}

\begin{remark} We remark that, like in the uncontrolled case, as $\Pu$ is a left stochastic matrix, then the total average on the graph only depends on the binary exchange process, i.e.,
\begin{equation}\label{eq:ave_tot_u}
    \dt \sum_{i\in \cI} \rho_i^u m_i^u= \mu \sum_{i\in \cI}(1-u_i^\mu) \rho_i^{u^2} m_i^u (\nu_2^i-\nu_1^i).
\end{equation}
The latter suggests the fact that a mere control on the interaction rate $\mu$ is sufficient in order to mitigate the evolution of the total weighted average (first moment). Moreover, as the microscopic processes are independent and the evolution of the physical quantity is mainly due to the binary exchange process, a mere control on the migration mechanism is surely not enough. Specifically, as already observed, it would not mitigate the total weighted average. However, as the traveling individuals carry on with them their state $v$ while migrating across the nodes, if the control was only exerted on the the interaction rate, then it would need to be stronger in order to mitigate the propagation of large values of $v$ across the graph. Conversely, these observations give a further motivation to the choice of exerting the two independent control mechanisms: in fact, controlling both processes with the same control could be over expensive.
\end{remark}

\subsection{Optimality conditions}\label{sec:optimality}
From now on, we drop the apex $^u$ on the average quantities $\rho_i, \rho_i m_i$ of the controlled problem~\eqref{eq:boltz.fi.c}.
We now want to find the optimal control $\bar{u}_i^\mu,\bar{u}_i^\chi$, by adapting a technique used in~\cite{albi2021}. Let us then discretize equation~\eqref{eq:boltz.fi.c} with a time step $h$
\(\label{eq:discreto}
\begin{aligned}
\int_\Rp \phi(v) f_i(t+h,v)\, dv &= \int_\Rp \phi(v) f_i(t,v)\, dv
                                  + h\chi \int_\Rp \phi(v) \bigl[ \sum_{j\in\cI}\Puij f_i(t,v) - f_i(t,v) \bigr] \, dv\\
                                 &{\hphantom{{}={}}}+ h \mu (1-u_i^\mu)\, \ev*{\int_\Rp\!\int_\Rp (\phi(v') - \phi(v)) f_i (t,v) f_i(t,v_*) \, dv_\ast \, dv}.
\end{aligned}
\)
We consider the cost functional:
\begin{equation}
\J_h^i(u_i^\chi,u_i^\mu, \rho_i m_i) = \psi(\rho_i m_i(t+h)) + \frac12\nu_i^\chi {u_i^\chi}^2+ \frac12\nu_i^\mu {u_i^\mu}^2,
\end{equation}
because we want to minimize the first moment $\rho_i m_i$, that is the average weighted by the mass in each node.
The minimization conditions are
\[
D_{u_i^\chi} \J_h^i(u_i^\chi,u_i^\mu, \rho_im_i) = 0,  \qquad D_{u_i^\mu} \J_h^i(u_i^\chi,u_i^\mu, \rho_i m_i) = 0, 
\]
that imply
\[
\frac{d \rho_i m_i}{d u_i^\chi} \psi'(\rho_i m_i(t+h)) + \nu_i^\chi u_i^\chi = 0, \qquad \frac{d \rho_i m_i}{d u_i^\mu} \psi'(\rho_i m_i(t+h)) + \nu_i^\mu u_i^\mu = 0.
\]
Exploiting equation~\eqref{eq:momentum_u}, the latter conditions  can be shown to be equivalent to
\[
\nu_i^\chi u_i^\chi + \psi'(\rho_i m_i(t+h)) h\left[-\chi \left(\sum_{\substack{j\in\cI,\\ j\neq i}} \Pij \rho_j m_j \right)\right] = 0, \quad \nu_i^\mu u_i^\mu+\psi'(\rho_i m_i(t+h)) h\left[-\mu (\nu_2^i-\nu_1^i)\rho_i^2 m_i\right] = 0.
\]
Now, if we impose $\nu_i^\alpha = h k_i^\alpha, \, \alpha=\chi, \mu$ for suitable $k_i^\chi, k_i^\mu>0$, we can write
\(
\bar{u}_i^\chi(t+h) = \psi'(\rho_i m_i(t+h)) \biggl[\frac{\chi}{k_i^\chi} \left(\sum_{\substack{j\in\cI,\\ j\neq i}} \Pij \rho_j m_j \right)\biggr], \qquad \bar{u}_i^\mu(t+h) = \psi'(\rho_i m_i(t+h)) \biggl[ \frac{\mu}{k_i^\mu}(\nu^i_2 - \nu^i_1)\rho_i^2 m_i \biggr].
\)

Letting $h \to 0$ we obtain that the \textit{optimal control problem} is defined by equation~\eqref{eq:boltz.fi.c} with 
\begin{equation}\label{eq:u}
u_i^\chi = \min\{\max\{\delta, \bar{u}_i^\chi\}, 1\}, \qquad u_i^\mu = \min\{\max\{0, \bar{u}_i^\mu\}, 1\},
\end{equation}
with $\delta > 0$ being the minimum entry of $\Pij$, and
\(\label{eq:u_star}
\bar{u}_i^\chi(t) = \psi'(\rho_i m_i(t)) \left[\dfrac{\chi}{k_i^\chi} \left(\sum_{\substack{j\in\cI,\\ j\neq i}} \Pij \rho_j m_j \right)\right],
\qquad \bar{u}_i^\mu(t) = \psi'(\rho_i m_i(t)) \left[\dfrac{\mu}{k_i^\mu} (\nu_2^i-\nu_1^i)\rho_i^2 m_i\right].
\)

In the same spirit as in~\cite{albi2021}, we consider 
\begin{equation}\label{def:psi}
\psi(x) = \dfrac{x^q}{q}, \quad q>1,
\end{equation}
which is $C^1$, so that the limit for $h\to 0^+$ in~\eqref{eq:u} is well defined.

Now, we discuss the compatibility of~\eqref{eq:u_star}. As the controls $u_i^\chi, u_i^\mu$ are actually defined by~\eqref{eq:u}, then they are compatible by construction, i.e.,
\[
u_i^\chi \in [0,1], \qquad u_i^\mu \in [0,1].
\]

However, we want to discuss possible ranges of values of $k_i^\mu, k_i^\chi$ allowing to obtain a compatible $\bar{u}_i^\chi, \bar{u}_i^\mu$, at least in the upper-bound. To this aim we impose $\bar{u}_i^\chi, \bar{u}_i^\mu\le 1$. Let us assume that $\psi'$ is positive. 
Differently with respect to~\cite{albi2021}, here the argument of $\psi$ is not monotone, but we determine a constraint on $k_i^\mu, k_i^\chi$ by imposing that $m_i(t)<m_i(0)$, as the final purpose is to diminish the average viral load, which is expected in the long run to be lower than the value it holds at time $t=0$. When the latter is not satisfied, the compatibility of the controls will be guaranteed by~\eqref{eq:u}. For the migration dynamics, the condition $\bar{u}_i^\chi\le 1$ can be satisfied when imposing $m_j(t)<m_j(0), \, \forall j$ by choosing for example
\begin{equation}\label{k_chi}
k_i^\chi \ge  \rho_i^{q-1 }(t)m_i^{q-1}(0) \chi  m_{\bar{i}}(0)(1-P_{ii}),  \qquad \bar{i}=\displaystyle \argmax_{j \in \cI-\lbrace i \rbrace} m_j(0),
\end{equation}
that can be obtained by remembering that $\rho_j(t) \le 1$ and using $1-\displaystyle\sum_{\substack{k\in \cI\\k\ne i}} P_{ki}=P_{ii}$ that follows from~\eqref{eq:P}.

We remark that when $\chi=0$ the penalization coefficient is not needed as there is no migration on the graph.
For the infection dynamics, the condition $\bar{u}_i^\mu\le 1$ when $m_i(t)<m_i(0)$ is satisfied by choosing for example
\begin{equation}\label{k_mu}
    k_i^\mu \ge \rho_i^{q+1}(t)\, m_i^{q-1}(t)m_i(0) \mu(\nu_2^i-\nu_1^i).
\end{equation}
The minimal choice determined by setting $k_i^\mu$ equal to the right hand side in the inequality~\eqref{k_mu} leads to
\(\label{eq:bar_u_conv}
\bar{u}_i^\mu=\dfrac{m_i(t)}{m_i(0)},
\)
that is smaller than one as long as $m_i(t)<m_i(0)$, while the interactions are stopped as soon as $m_i(t) \ge m_i(0)$.

We remark that also a larger choice for $k_i^\mu$ is possible, such as for example $k_i^\mu = \rho_i^{q+1}(t)m_i^q(0) \mu(\nu_2^i-\nu_1^i)$, which implies $\bar{u}_i^\mu=\left(\dfrac{m_i(t)}{m_i(0)}\right)^q$, or $k_i^\mu = m_i^q(0) \mu(\nu_2^i-\nu_1^i)$. These two choices automatically lead to a larger control on the interactions.
Moreover, if $\nu_2^i > \nu_1^i$, then the control $\bar{u}_i^\mu$ is positive, as well as the penalization coefficient; else if $\nu_2^i \le \nu_1^i$, then $\bar{u}_i^\mu \le 0$ so that $u_i^\mu =0$. In fact, if $\nu_2^i \le \nu_1^i$ the control is not needed, and coherently there is no actual constraint on $k_i^\mu$, and we have that $u_i^\mu =0$.

We remark that, if $\psi$ is a function of $m_i$, i.e., we control the average and not the weighted one $\rho_i m_i$, then we find, instead of $\bar{u}_i^\chi, \bar{u}_i^\mu$, other controls, that we shall denote $\tilde{u}_i^\chi,\tilde{u}_i^\mu$, that are given by
\(\label{eq:u_star_bis}
\tilde{u}_i^\chi = \psi'(m_i(t)) \left[\frac{\chi}{k_i^\chi} \left(\sum_{\substack{j\in\cI,\\ j\neq i}} \Pij \rho_j m_j \right)\right], \qquad \tilde{u}_i^\mu=\psi'(m_i(t))\left[\dfrac{\mu}{k_i^\mu} (\nu_2^i-\nu_1^i)\rho_i m_i\right].
\)
Then we have that
\[
\dfrac{\bar{u}_i^\alpha}{\tilde{u}_i^\alpha} = \dfrac{\psi'(\rho_i m_i) \rho_i}{\psi'(m_i)}, \quad \alpha=\chi, \mu.
\]
Given the choice~\eqref{def:psi}, then the latter ratio is exactly $\rho_i^q$. This is smaller than 1 and equates 1 only when $\rho_i=1$, i.e., when all the population is in node $i$. This is coherent with the fact that when all the population is in node $i$, controlling the population (and its average) on the entire graph and controlling only node $i$ is the same; at the same time, small values of the density $\rho_i$ in $i$ imply the fact that the value $m_i$ is not a reliable average as there are not many agents in $i$. These considerations justify our choice of exerting the control on $\rho_i m_i$ and not on $m_i$.

\paragraph{Global control}\label{sec:global}

Before turning to the analysis of the aggregate quantities of the optimal control model, we now justify the reason why it is more convenient to exert the proposed intra-node control instead of a global uniform control on the whole graph. When defining a global control, we aim at minimizing the the total weighted average (first moment), which is defined by
\[
\rho m(t) = \displaystyle\sum_{i\in\cI} \rho_i(t)m_i(t),
\]
and $\psi$ will depend on $\rho m$.
In particular, as a proof of concept, as the binary interactions are the stronger effect when it comes to the increase of the infection (even not to the diffusion), then we only consider the control on the binary interactions. Moreover, we remark that the global mean $\rho m$ is not affected by the mobility, as~\eqref{eq:ave_tot} holds.
Hence, if we consider a global multiplicative control on the interactions $\mu \to \mu\cdot (1-u)$, independent on specific nodes, we have
\(\label{eq:totalmean}
\dt \rho m(t) = \dt \sum_{i\in \cI} \rho_i(t)m_i(t) = \mu (1-u) \sum_{i\in \cI} (\nu_2^i - \nu_1^i)\rho_i(t)m_i(t), 
\)
where the control $u$ satisfies
\[
u = \psi'(\rho m) \frac{\mu}{k} \sum_{j\in \cI} (\nu_2^j - \nu_1^j)\rho_j(t)m_j(t),
\]
where, as mentioned before, $\psi$ depends on $\rho m$.
Conversely, if we consider a targeted intra-node control on the $i$-the node, (i.e., we impose $\mu \to \mu \cdot (1-\tilde u_i$)), still with the aim of minimizing $\rho m$, we obtain
\(\label{eq:totalmeantilde}
\dt \rho m(t) = \mu \sum_{i\in \cI} (1-\tilde u_i) (\nu_2^i - \nu_1^i)\rho_i(t)m_i(t), 
\)
where the following relation for $\tilde u_i$ needs to hold:
\[
\tilde u_i  = \psi'(\rho m) \frac{\mu}{k}  (\nu_2^i - \nu_1^i)\rho_i(t)m_i(t).
\]

We now remark that we have
\[
u=\sum_{i\in \cI} \tilde{u}_i,
\]
which implies $\tilde u_i < u$. Therefore, controlling the global weighted average when trying to minimize it with a global control, implies the fact that in some nodes the control may be too high. Instead, a localized control in each node is still efficient and has the effect of not penalizing everyone when it is not needed. This justifies the choice of implementing the control in each node. 

\subsection{The aggregate trend under optimal control}

The asymptotic behavior of the masses and of the averages heavily relies on the controls. While the trend of the averages will be investigated in the following, we now focus on the asymptotic trend of the masses. We highlight the fact that the evolution equation for the controlled masses may be written in vector notation as
\begin{equation}\label{eq:brho_u}
\frac{d}{dt} \brho^u = \chi [\Pu(t) - \cId]\brho^u,
\end{equation}
where $\Pu(t)$ is time dependent regardless of the specific choice of the optimality condition for the control on the mobility.
As a consequence, the existence and stability of the asymptotic stationary state of~\eqref{eq:brho_u} is not evident and cannot be investigated invoking the Perron-Frobenius theory.

One can, instead, exploit the theory of linear time variant systems with Metzler matrices~\cite{Metzler}. We recall that a Metzler matrix is a square time dependent matrix in which all the off diagonal entries are nonnegative after a given time.
Remarking that $\P$ and $\Pu$ are  Metzler matrices, one can prove the following result, that relies on assumptions on the basic mobility matrix $\P$.
\begin{prop}\label{prop:metzler}
    Let us suppose that the left stochastic transition matrix $\P$ satisfies
    \begin{equation}\label{hp_P}
        P_{ii} \ge \sum_{\substack{k \in \cI \\ k\ne i }}P_{ik} -1, \qquad \forall i \in \cI.
    \end{equation}
Therefore the system~\eqref{eq:brho_u} is globally asymptotically stable.
\end{prop}
\begin{proof}
The present proof relies on verifying the hypothesis of Theorem 4.1 in~\cite{Metzler}.
Let us define the Metzler matrix
\[
\Delta^M_{i,j}(t)=
\begin{cases}
    2(\Puij-1)+ \displaystyle\sum_{\substack{k \in \cI, \\ k\ne i }} P_{ik}^u \quad i= j,\\
    \Puij \quad i\ne j.
\end{cases}
\]
There exists a constant Metzler matrix $\Delta^{M_+}_{i,j}$ which is defined by
\[
\Delta^{M_+}_{i,j}(t)=
\begin{cases}
    -2+\displaystyle\sum_{\substack{k \in \cI, \\ k\ne i }} P_{ik} \quad i= j,\\
    P_{ij} \quad i\ne j.
\end{cases}
\]
that satisfies
\[
\Delta^M_{i,j}(t) \le \Delta^{M_+}_{i,j}, \, \forall i,j \in \cI, \quad \forall t> 0,
\]
because of the definition of $\Pu$ and exploiting $-u_i^\chi \le 0$.
Moreover, if~\eqref{hp_P} holds, then $\Delta^{M_+}$ has negative eigenvalues. As a consequence both the hypotheses of Theorem 4.1 in~\cite{Metzler} hold and we can conclude.
\end{proof}
\begin{remark}
    If $\P$ is right stochastic, then the hypothesis~\eqref{hp_P} is satisfied. 
\end{remark}
We now turn our attention to the existence of an equilibrium for the averages, which are the solutions of equations~\eqref{eq:momentum_u} and~\eqref{eq:mean_u} with~\eqref{eq:u},\eqref{eq:u_star},\eqref{def:psi},\eqref{k_chi},\eqref{k_mu}. We can consider the following three main cases.
\begin{itemize}
    \item If $\nu_1^i > \nu_2^i$ for all $i \in \cI$, then we immediately have that $m^u_i \to 0$ for all $i$. Indeed, since $u_i^\mu=0$, we have the monotonic decrease of the global mean
    \[
    \dt \sum_{\substack{i\in \cI}} \rho_i^u m_i^u = \sum_{i\in \cI} \mu(1 - u_i^\mu) (\nu_2^i - \nu_1^i) \rho_i^{u^2} m^u_i=\sum_{i\in \cI} \mu (\nu_2^i - \nu_1^i) \rho_i^{u^2} m^u_i < 0,
    \]
    whence we deduce that $\rho^u_i m^u_i \to 0$ for all $i \in \cI$, which gives the claim.
    \item If $\nu_1^i = \nu_2^i$ for all $i \in \cI$, then we may invoke Proposition~\ref{prop:metzler} replacing $\rho^u_i$ with the associate $\rho^u_i m^u_i$ to establish the existence of an asymptotically stable equilibrium point.
    \item If $\nu_1^i < \nu_2^i$ for all $i \in \cI$, then the global mean is monotonically increasing, and some more hypotheses are needed in order to also prove the existence and stability of equilibria. Proposition~\ref{prop:bounded} is devoted to this particular case. 
\end{itemize}

In the following proposition we drop the ${}^u$ apex and time dependence for convenience (except for the matrix~$\Pu$).
\begin{prop}\label{prop:bounded}
Let us consider $\nu_2^i > \nu_1^i$ for all $i \in \cI$ 
 and the specific choice for the control given by~\eqref{eq:bar_u_conv}. Let us assume the following
\begin{enumerate}[(i)]
\item\label{cond:mi} For all $i \in \cI$ there exists a time $t_2^i > 0$ such that, for all $t \ge t_2^i$, holds
\(
- m_i(t_2^i) <
        \chi \int_{t_2^i}^t \biggl(\sum_{\substack{j\in \cI}} \Puij(s) \frac{\rho_j(s)}{\rho_i(s)} ( m_j(s)  - m_i(s)) \biggr)\, ds.
\)
\item\label{cond:rhoimi} For all $i \in \cI$ there exists a constant $\alpha_i > 1$, a time $t_1^i > 0$ and a function $r(\; \cdot\;)$ such that, for all $t \ge t_1^i$, holds
\begin{gather}
    \label{eq:r-infinito}
    r(t) \ge 0 \text{ and } \lim_{t \to +\infty} r(t) = + \infty,\\
    \label{eq:bazooka-rhoimi}
    -\rho_i m_i(t_1^i) <
        \chi \int_{t_1^i}^t \biggl(\sum_{\substack{j\in \cI}} \Puij(s) \rho_j(s) m_j(s)  - \rho_i(s) m_i(s) \biggr)\, ds <
        \rho_i m_i(t_1^i) \Bigl(\frac{1}{r(t)} + \alpha_i - 1 \Bigr).
\end{gather}
\end{enumerate}

Then, we have that $\rho_i m_i$ and $m_i$ solutions to~\eqref{eq:momentum_u} and~\eqref{eq:mean_u} with~\eqref{eq:u}-\eqref{eq:u_star}-\eqref{def:psi}-\eqref{k_chi}-\eqref{eq:bar_u_conv} converge, eventually in a monotonic fashion, to finite equilibria $\rho_i^\infty m_i^\infty$ and $m_i^\infty$ for all $i \in \cI$.
\end{prop}
\begin{proof}
For the readability of the paper, the Proof is reported in Appendix~\ref{appendix:proof}.
\end{proof}
The results on the trend for masses and averages reported in Propositions~\ref{prop:metzler} and~\ref{prop:bounded} imply that Proposition~\ref{prop2} may be rephrased as
\begin{prop}\label{prop2_bis}
Assume the graph is strongly connected and that $\nu_1^i=\nu_1, \, \nu_2^i=\nu_2$. When $t\rightarrow \infty$, the stationary state $m_i$ of~\eqref{eq:mean_u}, where the controls are defined by~\eqref{eq:u}-\eqref{eq:u_star} with~\eqref{def:psi}-\eqref{k_chi}-\eqref{k_mu}, satisfies one of the following cases.
\begin{enumerate}
    \item If $\nu_1 > \nu_2$ then $m_i\rightarrow 0$ that is stable and attractive  for all $i \in \cI$;
     \item If $\nu_1 = \nu_2$ and under the hypotheses of Proposition~\ref{prop:metzler}, then $m_i \rightarrow m_i^\infty$ (to be determined) that is stable  for all $i \in \cI$;
    \item If $\nu_1 < \nu_2$ and under the hypotheses of Proposition~\ref{prop:bounded}, then $m_i \rightarrow m_i^\infty$ (to be determined) that is stable  for all $i \in \cI$.
   \end{enumerate}
\end{prop}
\begin{remark}
In the case $\nu_1^i = \nu_2^i$ for all $i \in \cI$, we may exchange the hypotheses of Proposition~\ref{prop:metzler} on the matrix $\P$ with a condition on $\rho_i m_i$ of the form~\ref{cond:rhoimi} (first inequality) and~\ref{cond:mi} of Proposition~\ref{prop:bounded}, and follow the proof of Proposition~\ref{prop:bounded} in order to prove the existence of asymptotically stable equilibrium points for the first order moments and averages.
\end{remark}

We now want to show the interplay of the intrinsic mitigating effect of the network (as shown in Proposition~\ref{prop2}) and of the control.  
We can rephrase Proposition~\ref{prop2} for the control problem as follows.
\begin{prop}\label{prop4bis} Let the graph be strongly connected and let us assume that $\mu=\dfrac{1}{\varepsilon}$, $\chi=1$, where $\varepsilon \ll 1$ in~\eqref{eq:boltz.fi.c}. Then, one of the following cases hold.
\begin{enumerate}
    \item Assume $\nu_1^i>\nu_2^i$ for all $i \in \cI-\lbrace i_\ast\rbrace$ and $\nu_1^{i_*}=\nu_2^{i_*}$. Then $u_i^\mu=u_i^\chi=0$ and $m_i \rightarrow 0$ for all $i \in \cI$, in particular also for $i=i_\ast$.
    \item Assume $\nu_1^i=\nu_2^i$ for all $i \in \cI-\lbrace i_\ast\rbrace$ and $\nu_1^{i_*}>\nu_2^{i_*}$. Let moreover $\cI_\ast=\lbrace j \in \cI: P_{i_\ast,j}>0\rbrace$. Then $m_i \rightarrow 0$ and $u_i^\chi \rightarrow 0$ for all $i \in \cI_*$.
    \item Assume $\nu_1^{i_*}<\nu_2^{i_*}$ for some $i_\ast \in \cI$. Then $m_i \rightarrow m_i^\infty<\infty$ when $t \rightarrow \infty$ for all $i \in \cI$ s.t. $\nu_1^i \le \nu_2^i$, while $m_i \rightarrow 0$ when $t\rightarrow \infty$ for all $i \in \cI$ s.t. $\nu_1^i>\nu_2^i$, so that $u_i^\chi \rightarrow 0$.
\end{enumerate}
\end{prop}
\begin{proof}
The proof follows straightforwardly from the corresponding points in Proposition~\ref{prop2}, and from the definition of~\eqref{eq:u}-\eqref{eq:u_star} with~\eqref{k_chi}-\eqref{k_mu}.  
\end{proof}
We can conclude that, as infection and healing are both consequence of the same microscopic process, that is a binary interaction, the control defined for mitigating the infection process also acts on the healing process, and this is counterproductive. As a consequence, in the next section we shall propose a new kinetic model separating the infection and healing processes.

\section{A kinetic model with infection and healing}\label{sec:healing}

 We now propose a new kinetic model in which the infection and the healing are modeled by two different independent microscopic processes. The infection will be a consequence of binary interactions, while the healing will be modeled as an autonomous linear process. 
\subsection{The kinetic model}
For the latter considerations, we consider the kinetic model defined by
\(
\begin{aligned}\label{eq:boltz.fi_new}
\frac{d}{dt}\int_{\Rp} \phi(v) f_i(t,v) \, dv
&= \chi \int_{\Rp} \phi(v) \biggl(\sum_{j \in \cI}\Pij f_j(v,t) -f_i(t,v)\biggr)\, dv\\
&+ \sigma \int_{\Rp} \phi(v)Q^+(f_i,f_i)(v,t)\, dv+ \gamma \int_{\Rp} \phi(v)Q^-(f_i)(v,t)\, dv, \qquad i\in \cI,
\end{aligned}
\)
where $\phi\colon \Rp \to \R$ is a test function and $Q^+(f_i,f_i),Q^-(f_i)$ are two different \emph{collision-like operators}, defined as
\(
\int_{\Rp} \phi(v) Q^+(f_i,f_i)(v,t)\, dv = \int_{\Rp}\!\int_{\Rp} \ev{\phi(v') - \phi(v)} f_i(t,v)f_i(t,v_*)\, dv_* \, dv, 
\)
where
\begin{equation}\label{v_prime_new}
 v'=v+\nu_2^i v_\ast +\eta_1' v, \quad v_\ast'= v_\ast+\nu_2^i v +\eta_2' v_\ast  ,
\end{equation}
while
\(
\int_{\Rp} \phi(v) Q^-(f_i)(v,t)\, dv = \int_{\Rp}\! \ev{\phi(v'') - \phi(v)} f_i(t,v)\, dv, 
\)
where
\begin{equation}\label{v_second}
 v''=v-\nu_1^i v +\eta'' v, 
\end{equation}
being $\eta_1', \eta_2', \eta''$ white noises.
We remark that $Q^+$ implements a binary interaction leading to infection and thus represents a \emph{gain} for the viral--load, while $Q^-$ is a linear process modeling the autonomous and spontaneous healing and is a \emph{loss} term for the viral--load. The parameters $\sigma$ and $\gamma$ are the corresponding frequencies.

\subsection{Aggregate trend}
The evolution of the masses~\eqref{eq:rho_i} is the same as for the model~\eqref{eq:boltz.fi} and the same result of Proposition~\ref{prop:rho^inf} holds.
Conversely, setting $\varphi(v)=v$ in~\eqref{eq:boltz.fi_new} we may investigate the evolution of the first moment $\rho_i m_i$ 
which turns out to satisfy the equation
\begin{equation}
	\frac{d}{dt}(\rho_im_i)=\chi\left(\sum_{j\in\cI}P_{ij}\rho_jm_j-\rho_im_i\right)+(\sigma\nu_2^i \rho_i-\gamma\nu_1^i)\rho_im_i.
	\label{eq:rho_im_i_new}
\end{equation}
Moreover, the evolution of the average 
$$ m_i(t)\coloneqq\frac{1}{\rho_i(t)}\int_{\R_+}vf_i(t,v)\,dv, $$
is 
\begin{equation}\label{eq:momentum_new}
	\dt m_i=\frac{\chi}{\rho_i}\sum_{j\in\cI}P_{ij}\rho_j(m_j-m_i)+(\sigma\nu^i_2\rho_i-\gamma \nu^i_1)m_i.
\end{equation}
As expected, the average varies because of the sum of three independent contributions: the migration on the graph, the infection (increase of viral load) and the healing (decrease of viral load), as
\begin{equation}
\dt m_i=\left[\dt m_i\right]_\chi +\left[\dt m_i\right]_\sigma+\left[\dt m_i\right]_\gamma,
\end{equation}
where $\left[\dt m_i\right]_\chi$ is defined in~\eqref{eq:ave_chi}, while
\begin{equation}\label{eq:ave_mu_new}
\left[\dt m_i\right]_\sigma=\sigma\nu^i_2\rho_im_i, \qquad \left[\dt m_i\right]_\gamma=-\gamma\nu^i_1 m_i
\end{equation}
are the contributions of the infection and of the healing processes, respectively.
In particular, we remark that, aside of the term $\left[\dt m_i\right]_\chi$, the increase of $m_i$ depends on the infection rate $\sigma \nu_2^i$ and, as expected, on the mass $\rho_i$ in the $i-$th node; conversely the decrease of $m_i$ only depends on the healing process and on the rate $\gamma \nu_1^i$.

\begin{figure}[htbp]
    \hbox to \textwidth{%
    \includegraphics[width = 0.45\linewidth]{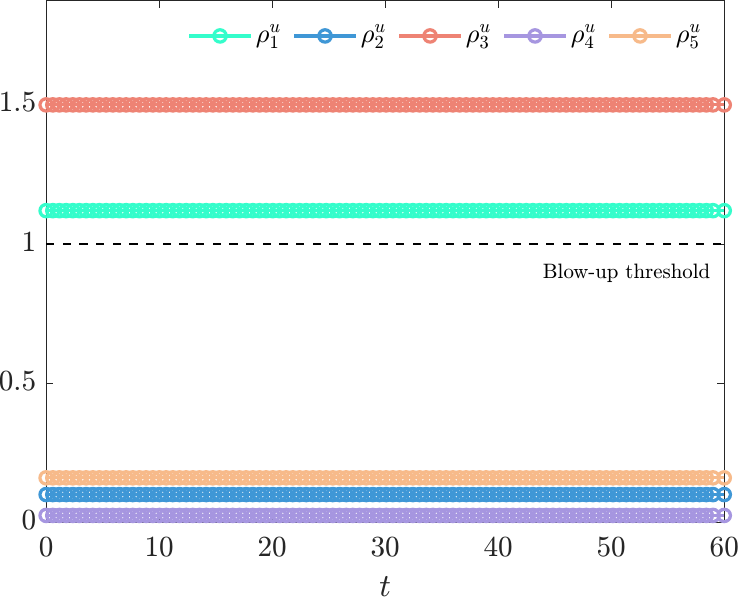}\hfil
    \includegraphics[width = 0.45\linewidth]{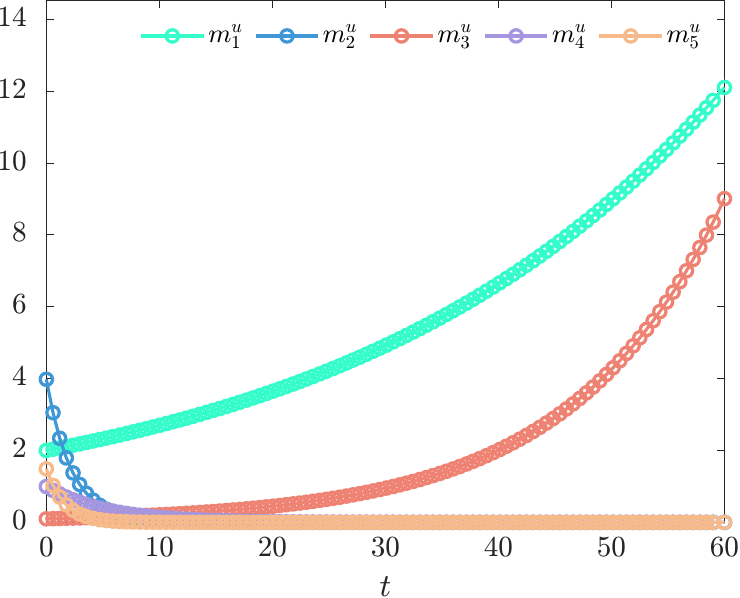}}
    \caption{Left to right: evolution in time of number of agents (divided by $\rho_i^c$) and average viral load, in absence of mobility ($\chi = 0$). We see that, in those nodes where the initial mass fraction  $\rho_i$ is greater than the associated critical value $\rho_i^c$, the average viral load grows exponentially, while in the other nodes vanishes. Refer to Section~\ref{sec:num2} for additional details about the simulation.}
    \label{fig:rhoic}
\end{figure}

If there is no graph ($n=1$) or there is no migration on the graph ($\chi=0$) we can remark that $\rho_i(t) \equiv\rho_i^\infty=\rho_i(0), \, \forall \, t>0, \, \forall i\in \cI$. Let us define
    \begin{equation}\label{def:rhoic}
        \rho_i^c=\dfrac{\gamma \nu_1^i}{\sigma \nu_2^i}.
    \end{equation} 
    We can distinguish two cases.
\begin{itemize}
    \item When $\rho_i^c\ge 1$, i.e., $\sigma \nu_2^i \le \gamma \nu_1^i$, then $m_i \rightarrow 0$ as, even if $\sigma \nu_2^i = \gamma \nu_1^i$, then $\sigma \nu_2^i \rho_i \le \gamma \nu_1^i$, being $\rho_i \le 1$.
    \item The case $\rho_i^c<1$, i.e., $\sigma \nu_2^i > \gamma \nu_1^i$, is more complex. According to the value of $\rho_i^\infty$ with respect to $\rho^c_i$, the evolution of the average may invert its trend naturally (without control) with respect to the initial one. In fact, if 
    $$
    (\rho_i^c-\rho_i(0))(\rho_i^c-\rho_i^\infty)<0, 
    $$
    and $\rho_i(0)>\rho_i^c$,
    then $\exists t_0 >0$ such that for $t<t_0$, then $\rho_i(t)>\rho_i^c$, while for $t>t_0$ $\rho_i(t)<\rho_i^c$. For $t<t_0$, we have that $\sigma\nu_2^i \rho_i-\gamma\nu_1^i >0$ and, then, $m_i$ increases; conversely, when $\rho_i<\rho_i^c$ for $t>t_0$, then $\sigma \nu_2^i \rho_i-\gamma\nu_1^i <0$ and $m_i$ decreases. See Figure~\ref{fig:rhoic} for differences in the time evolution of $m_i(t)$ depending on the ratio $\rho_i^\infty/\rho_i^c$.
\end{itemize}

On the other hand, if there is migration on the graph, we anyway have that, as a consequence of~\eqref{eq:P}, the variation of the total average on the graph is only due to the interactions inside the nodes, as
\begin{equation}\label{eq:ave_lin_new}
   \dt \sum_{i\in \cI} \rho_i m_i= \sigma \sum_{i\in \cI} \rho_i m_i (\sigma\nu_2^i\rho_i - \gamma\nu_1^i).
\end{equation} 

As a consequence, again, it is possible to prove by means of a linear stability analysis that
\begin{prop}\label{prop11}
Assume the graph is strongly connected. When $t\rightarrow \infty$, one of the following statements hold true for the solution $m_i$ of~\eqref{eq:momentum_new}.
\begin{enumerate}
    \item If $\rho_i^\infty < \rho_i^c$ $\forall i$ then $m_i\rightarrow 0$   $\forall i \in \cI$.
    \item If $\gamma \nu_1^i>\sigma \nu_2^i$, i.e., $\rho_i^c>1$, then $m_i\rightarrow 0$  $\forall i \in \cI$.
    \item If $\rho_i^c < \rho_m^\infty\coloneqq \displaystyle\min_{i\in\cI} \rho_i^\infty$ $ \forall i\in \cI$ then $m_i \rightarrow \infty$ for some $i \in \cI$.
\end{enumerate}
\end{prop}
\begin{proof}
We are interested
in the stability of the asymptotic state $\bmm^\infty=0$ which represents the eradication of the infection
in all nodes of the network. Then we now consider the linearization of equation~\eqref{eq:ave_lin_new} around the equilibrium configuration $(\brho,\boldsymbol{m})=(\brho^\infty, \boldsymbol{0})$. Writing
$$ \rho_i=\rho^\infty_i+\epsilon\tilde{\rho}_i, \qquad m_i=\epsilon\tilde{m}_i, $$
where $\epsilon>0$ is a small parameter, and plugging it into~\eqref{eq:ave_lin_new}, we obtain
\begin{equation}\label{eq:ave_lin}
   \dt \sum_{i\in \cI} \rho_i^\infty \tilde{m}_i=  \sum_{i\in \cI} \rho_i^\infty \tilde{m}_i (\sigma\nu_2^i\rho_i^\infty-\gamma\nu_1^i).
\end{equation}  
If $\rho_i^\infty < \rho_i^c \, \forall i \in \cI$, then setting $d\coloneqq \displaystyle \max_{i\in\cI}\left( \sigma\nu_2^i\rho_i^\infty-\gamma\nu_1^i\right)$, we have that
\[
\dt \sum_{i\in \cI} \rho_i^\infty \tilde{m}_i\le d \sum_{i\in \cI} \rho_i^\infty \tilde{m}_i, 
\]
and $m_i\rightarrow 0$ as $d < 0$. This proves the point 1 of the Proposition.
Moreover, as $\rho_i^\infty \le 1$, we have that 
\(
\dt \sum_{i\in \cI} \rho_i^\infty \tilde{m}_i\le  \sum_{i\in \cI} \rho_i^\infty \tilde{m}_i (\sigma\nu_2^i-\gamma\nu_1^i).
\)
Then, if $\sigma\nu_2^i-\gamma\nu_1^i<0 \, \forall i\in \cI$, we can conclude again that $m_i\rightarrow 0$ $\forall i \in \cI$. This proves point 2.
Conversely, if $\rho_i^c < \rho_m^\infty \coloneqq \displaystyle \min_{i\in\cI} \rho_i^\infty$, then
\(
\dt \sum_{i\in \cI} \rho_i^\infty \tilde{m}_i\ge  \sum_{i\in \cI} \rho_i^\infty \tilde{m}_i b,
\)
where $b\coloneqq\displaystyle \min_{i\in\cI} \left(\sigma\nu_2^i \rho_m^\infty-\gamma\nu_1^i\right)$ that is a positive quantity. In conclusion, $\exists i \in \cI $ such that $m_i\rightarrow \infty$.
\end{proof}

\begin{remark}
    Notice that, if $\rho_i^c>1$ then $m_i \rightarrow 0$ (case 2). If $\rho_i^c <1$ then if Case 1 holds then there is eradication, while if Case 3 holds then there is blow up. 
\end{remark}

\begin{remark}
    We do not consider the case of constant infection/healing coefficients, as now, also in the linearized case, the infection coefficient always depends on the node through $\rho_i^\infty$.
\end{remark}

\begin{remark}
    Notice that this dynamics is influenced by the values of $\rho_i^\infty, \rho_i^c, \sigma \nu_2^i, \gamma \nu_1^i$, but, also, by $\rho_i(0)$ and by the migration strategy defined by $\P$ and $\chi$. For this reason, the values of the averages $m_i$ may also oscillate in some parameters regimes.
Specifically, we may consider the following cases:
\begin{itemize}
    \item when $\sigma \nu_2^i \le \gamma \nu_1^i \, \forall i \in \cI$, then $m_i \rightarrow 0$ $\forall i \in \cI$, as this corresponds to $\rho_i^c \ge 1 \ge \rho_i^\infty, \, \forall i \in \cI$;
    \item when $\sigma \nu_2^i > \gamma \nu_1^i$ for some $i$, then we have that, according to the value of $\rho_i(t)$ with respect to $\rho^c_i$, the evolution of the average may invert its trend. In fact, when $\rho_i(t)>\rho_i^c$, then $\sigma\nu_2^i \rho_i(t)
    -\gamma\nu_1^i >0$, and then the contribution of the exchange process to $m_i$ ($\left[\dt m_i\right]_\sigma+\left[\dt m_i\right]_\gamma$) is positive, while when $\rho_i(t)<\rho_i^c$, then $\sigma\nu_2^i \rho_i(t)-\gamma\nu_1^i <0$, and this contribution is negative. In particular, we can argue that, given an initial condition $\rho_i(0)$ in each node and the corresponding stationary state $\rho_i^\infty$, then, if
    \[
    (\rho_i(0)-\rho_i^c)(\rho_i^\infty-\rho_i^c)<0
    \]
    then there exists $\bar{t}$ such that $\rho_i(\bar{t})=\rho_i^c$. Then, if $\rho_i(0)>\rho_i^c$, then $m_i \rightarrow 0$, else if $\rho_i(0) < \rho_i^c$, then it might happen that $m_i \rightarrow \infty$.
    Conversely, when
    \[
    (\rho_i(0)-\rho_i^c)(\rho_i^\infty-\rho_i^c)>0
    \]
    then, if $\rho_i^\infty, \rho_i(0) <\rho_i^c$ then $m_i\rightarrow 0$. Else, if $\rho_i^\infty, \rho_i(0) >\rho_i^c$, then it might happen that $m_i\rightarrow \infty$.
\end{itemize}
\end{remark}

\section{Optimal control of the kinetic model with healing}\label{sec:healing.control}
In this section we implement the optimal control on the kinetic model with healing. In particular, we include a control on the mobility, as in Section~\ref{sec:control}, and on the infection process ruled by binary interactions, while the autonomous healing process will not be controlled.
\subsection{Control in each node}
We now implement two independent controls on the microscopic migration dynamics and on the microscopic binary interactions leading to infection. 
We, then, define the controlled problem as
\begin{equation}\label{eq:boltz.fi.c_new}
\begin{aligned}
\dt \int_\Rp \phi(v) f_i(t,v)\, dv &= \chi \int_\Rp \phi(v) \left[ \sum_{j \in \cI}\Puij f_j(t,v) - f_i(t,v) \right] \, dv\\
                                          &\hphantom{{}=}+ \sigma_i \int_{\Rp} \phi(v)Q^+(f_i,f_i)(v,t)\, dv+ \gamma \int_{\Rp} \phi(v)Q^-(f_i)(v,t)\, dv,
\end{aligned}
\end{equation}
where the control matrix $\Puij$ is defined by equation~\eqref{eq:defpuij}, while the control $\sigma_i$ on the binary interactions is defined by
\begin{equation}\label{def:sigmai}
    \sigma_i\coloneqq\sigma (1-u_i^\sigma)
\end{equation}
and has the effect of reducing the infection rate inside each node. 

 Concerning the controlled macroscopic equations, we have that the evolution of the masses is again given by~\eqref{eq:brho_u}, while
setting $\varphi=v$ in~\eqref{eq:boltz.fi.c_new} we obtain for the controlled weighted average $\rho_i^u m_i^u$
\(\label{eq:momentum_u_new}
\dt \rho_i^u(t)m_i^u(t) = \chi \biggl[\sum_{j\in \cI} \Puij \rho_j^u(t)m_j^u(t) - \rho_i^u(t)m_i^u(t) \biggr]
                      + \sigma (1-u_i^\sigma) \nu^i_2\rho_i^{u^2}(t)m_i^u(t) - \gamma\nu^i_1\rho_i^{u}(t)m_i^u(t),
\)
and, also
\(\label{eq:mean_u_new}
\begin{aligned}
\dt m_i^u(t)  = \chi \biggl[\sum_{\substack{j\in\cI \\ j\neq i}} \Puij \frac{\rho_j^u(t)}{\rho_i^u(t)}\bigl(m_j^u(t) - m_i^u(t)\bigr) \biggr]
                      + \sigma (1-u_i^\sigma) \nu^i_2 \rho_i^u(t)m_i^u(t) - \gamma \nu^i_1 m_i^u(t),
\end{aligned}
\)
i.e.
\[
\dt m_i^u=(1-u_i^\chi) \left[\dt m_i\right]_\chi + (1-u_i^\sigma) \left[\dt m_i\right]_\sigma-\left[\dt m_i\right]_\gamma ,
\]
where in the right hand side the quantities are defined in~\eqref{eq:ave_chi}-\eqref{eq:ave_mu_new}.
\begin{remark} Again, we can remark that for both the evolution of the masses and of the averages, the controls $1-u_i^\chi, 1-u_i^\sigma \in [0,1]$ reduce the time variation rate, without inverting the natural trend, with $u_i^\sigma,u_i^\chi=1$ stopping the time evolution. The difference here is that when increasing $u_i^\sigma$, the effect of $+\left[\dt m_i\right]_\gamma$ is stronger than the infection process, and it is not reduced like in the previous model of Section~\ref{sec:control}.
\end{remark}

\subsection{Optimality conditions}
From now on, we drop the apex $^u$ on the average quantities $\rho_i, \rho_i m_i$ of the controlled problem.
We now want to find the optimal control $\bar{u}_i^\sigma,\bar{u}_i^\chi$, as done in Section~\ref{sec:control}. Let us then consider a discretization in time of~\eqref{eq:boltz.fi.c_new} (being $h$ the time step)
\(\label{eq:discreto_new}
\begin{aligned}
\int_\Rp \phi(v) f_i(t+h,v)\, dv &= \int_\Rp \phi(v) f_i(t,v)\, dv
                                  + h\chi \int_\Rp \phi(v) \left[\sum_{j\in \cI} \Puij f_i(t,v) - f_i(t,v) \right] \, dv\\
                                 &{\hphantom{{}={}}}+ h \sigma (1-u_i^\sigma)\, \ev*{\int_\Rp\!\int_\Rp (\phi(v') - \phi(v)) f_i (t,v) f_i(t,v_*) \, dv_* \, dv}\\
                                 &{\hphantom{{}={}}}- h \gamma \int_{\Rp}\! \ev{\phi(v'') - \phi(v)} f_i(t,v)\, dv.
\end{aligned}
\)

We consider again the cost functional aiming at minimizing the average weighted by the mass in each node $\rho_i m_i$
\begin{equation}
\J_h^i(u_i^\chi,u_i^\sigma, \rho_i m_i) = \psi(\rho_i m_i(t+h)) + \frac12\nu_i^\chi {u_i^\chi}^2+ \frac12\nu_i^\sigma {u_i^\sigma}^2.
\end{equation}

The minimization conditions are
\[
D_{u_i^\chi} \J_h^i(u_i^\chi,u_i^\sigma, \rho_i m_i) = 0,  \qquad D_{u_i^\sigma} \J_h^i(u_i^\chi,u_i^\sigma, \rho_i m_i) = 0, 
\]
that imply
\[
\frac{d \rho_i m_i}{d u_i^\chi} \psi'(\rho_i m_i(t+h)) + \nu_i^\chi u_i^\chi = 0, \qquad \frac{d \rho_i m_i}{d u_i^\sigma} \psi'(\rho_i m_i(t+h)) + \nu_i^\sigma u_i^\sigma = 0,
\]
that, from equation~\eqref{eq:mean_u_new}, are equivalent to
\[
\nu_i^\chi u_i^\chi + \psi'(\rho_i m_i(t+h)) h\left[-\chi \left(\sum_{\substack{j\in\cI,\\ j\neq i}} \Pij \rho_j m_j \right)\right] = 0, \qquad \nu_i^\sigma u_i^\sigma+\psi'(\rho_i m_i(t+h)) h\left[-\sigma \nu_2^i \rho_i^2 m_i\right] = 0.
\]
Now, if we impose $\nu_i^\alpha = h k_i^\alpha, \, \alpha=\chi, \sigma$ for suitable $k_i^\chi, k_i^\sigma>0$, we can write
\(
\bar{u}_i^\chi(t+h) = \psi'(\rho_i m_i(t+h)) \biggl[\frac{\chi}{k_i^\chi} \left(\sum_{\substack{j\in\cI,\\ j\neq i}} \Pij \rho_j m_j \right)\biggr], \qquad \bar{u}_i^\sigma(t+h) = \psi'(\rho_i m_i(t+h)) \biggl[ \frac{\sigma}{k_i^\sigma}\nu^i_2 \rho_i^2 m_i \biggr].
\)
Then, when $h \to 0$ in equation~\eqref{eq:discreto_new}, we obtain that the optimal control problem now is defined by~\eqref{eq:boltz.fi.c_new} with 
\begin{equation}\label{eq:u_new}
u_i^\chi = \min\{\max\{\delta, \bar{u}_i^\chi\}, 1\}, \qquad u_i^\sigma = \min\{\max\{0, \bar{u}_i^\sigma\}, 1\},
\end{equation}
with $\delta > 0$ being the minimum entry of $\Pij$ and
\(\label{eq:u_star_new}
\bar{u}_i^\chi(t) = \psi'(\rho_i m_i(t)) \left[\dfrac{\chi}{k_i^\chi} \left(\sum_{\substack{j\in\cI,\\ j\neq i}} \Pij \rho_j m_j \right)\right], \qquad \bar{u}_i^\sigma(t) = \psi'(\rho_i m_i(t)) \left[\dfrac{\sigma}{k_i^\sigma} \nu_2^i\rho_i^2 m_i\right].
\)
Now, we discuss the compatibility of $\bar{u}_i^\sigma$, considering, again $\psi$ given by~\eqref{def:psi}.  Then, imposing that $\bar{u}_i^\sigma \le 1$, we find that a suitable choice is
\begin{equation}\label{eq:k_q+1}
k_i^\sigma \ge \rho_i^{q + 1}(t)m_i^q(t) \nu_2^i \sigma.
\end{equation}
However, now, we also impose that
\[
\sigma(1-u_i^\sigma)\nu_2^i <\gamma \nu_1^i,
\]
because in this way we choose a control that is strong enough in order to make the binary interaction process weaker than the healing one in each node.
The latter is satisfied if
\begin{equation}\label{cond_int}
\bar{u}_i^\sigma \ge 1-\rho_i^c.
\end{equation}
Of course, ~\eqref{cond_int} is meaningful only in the case $\rho_i^c <1$, which is the dangerous one. 
We have that~\eqref{cond_int} is satisfied when choosing
\begin{equation}\label{eq:comp_new_max}
 k_i^\sigma \le \frac{\rho_i^{q + 1}(t) m_i(t)^q \nu_2^i\sigma}{1 - \rho_i^c}.
\end{equation}
Therefore
\begin{equation}\label{k_interval}
k_i^\sigma(t) \in \biggl[
\rho_i^{q + 1}(t) m_i^q(t) \nu_2^i\sigma,
\frac{\rho_i^{q + 1}(t) m_i^q(t) \nu_2^i\sigma}{1 - \rho_i^c}
\biggr],
\end{equation}
that is well defined as $\rho_i^c < 1$ is required by~\eqref{cond_int}.
\begin{remark}
    Differently with respect to Section~\ref{sec:optimality}, we require~\eqref{eq:k_q+1} instead of 
    $k_i^\sigma \ge \rho_i^{q+1}(t)^qm_i(0) \nu_2^i \sigma$. This choice is linked to the necessity of having a good definition of the interval~\eqref{k_interval} for $k_i^\sigma$. This also leads to the fact that the penalization coefficient now depends on time also through the average. Specifically for the choice $k_i^\sigma=\frac{\rho_i^{q + 1}(t) m_i^q(t) \nu_2^i\sigma}{1 - \rho_i^c}$, i.e., the upper bound of the interval, this is reasonable. In fact, if the trend of the average $m_i$ is decreasing then restrictions can be made lighter: this corresponds to having a lower penalization coefficient in such a way that a too high control is not applied.

\end{remark}

\subsection{Aggregate trend of the optimal control problem}
The macroscopic quantities then evolve as~\eqref{eq:brho_u}-\eqref{eq:momentum_u_new} (or~\eqref{eq:mean_u_new}) with the controls defined by~\eqref{eq:u_new}-\eqref{eq:u_star_new} with~\eqref{k_chi}-\eqref{k_interval}. Again, in order to study the (at least linear) stability of possible stationary states $(\brho^\infty, \boldsymbol{0})$ we need to analyze the possibility of having a stationary transition probability. The latter amounts to investigating if the control on the mobility reaches a stationary state, i.e.
\[
\lim_{t\rightarrow \infty }u_i^\chi(t) = u_i^{\chi,\infty},
\]
 that is verified if $\rho_i m_i$ tends in time to a stationary finite value $\rho_i^\infty m_i^\infty$.
We consider the following possible cases.
\begin{itemize}
 \item When $\sigma\nu_2^i\le \gamma\nu_1^i$, then we know that $\rho_i^\infty m_i^\infty =0$ $\forall i \in \cI$.
\item When $\sigma\nu_2^i>\gamma\nu_1^i$ then $u_i^\sigma$ is chosen, thanks to~\eqref{k_interval}, in such a way  that $\rho_i^\infty m_i^\infty = 0$ $\forall i \in \cI$. 
\end{itemize}

In conclusion, the quantity $u_i^{\chi,\infty}$ defines a stationary controlled transition matrix $\P^{u,\infty}$ that is also irreducible. Again, it is possible to have a stationary state $\brho^\infty$ and Propositions~\ref{prop:rho^inf} holds true. We can also prove the controlled version of Proposition~\ref{prop11}.
\begin{prop}\label{prop12}
Assume the graph is strongly connected. When $t\rightarrow \infty$ the solution of~\eqref{eq:momentum_u_new} with~\eqref{eq:u_new}-\eqref{eq:u_star_new}-\eqref{k_chi}-\eqref{eq:comp_new_max} satisfies
$$
m_i \rightarrow 0,
$$
for any value of $\rho_i^c$.
\end{prop}
\begin{proof}
If $\rho_i^c \ge 1$ for a given $i\in\cI$, then we already knew that for the non-controlled problem $m_i\rightarrow 0$. If $\rho_i^c < 1$ for a given $i\in\cI$ we have that then $m_i \rightarrow 0$ thanks to the choice of the control defined by~\eqref{eq:u_star_new}-\eqref{k_chi}-\eqref{eq:comp_new_max}. 
\end{proof}

\begin{remark}
    The choice of controlling and minimizing the weighted average $\rho_i m_i$ inside each node is again due to the fact that the observations of Section~\ref{sec:optimality} and Section~\ref{sec:global} can be shown to hold true also in the present case.
\end{remark}

\subsection{Basic reproduction number}
In the same spirit as in~\cite{loytosin2021}, we can determine a basic reproduction number $\cR_0\geq 0$. The latter is a quantity that is typically defined in compartmental models and that represents the mean number of secondary infections caused by a single infected individual in a population of susceptible individuals. Mathematically, this parameter is defined as the one discriminating between the (linear) stability (if $\cR_0<1$) or instability (if $\cR_0>1$) of the disease-free equilibrium (no infected individuals).

As done in~\cite{loytosin2021}, as we are dealing with a non-compartmental viral load-based model, we define $\cR_0$ exploiting the stability/instability of the asymptotic state $\bmm^\infty=\bzero$. We then linearize~\eqref{eq:momentum_new} around the equilibrium $(\brho,\,\bmm)=(\brho^\infty,\,\bzero)$. Writing
$$ \rho_i=\rho^\infty_i+\epsilon\tilde{\rho}_i, \qquad m_i=\epsilon\tilde{m}_i, $$
where $\epsilon>0$ is a small parameter, and plugging into~\eqref{eq:momentum_new} we obtain, at the leading order in $\epsilon$, the following system of equations for the perturbations $\tilde{m}_i$:
$$ \frac{d\tilde{m}_i}{dt}=\frac{\chi}{\rho^\infty_i}\sum_{j\in\cI}P_{ij}\rho^\infty_j(\tilde{m}_j-\tilde{m}_i)+\mu(\nu_2\rho^\infty_i-\nu_1)\tilde{m}_i, \qquad i\in\cI. $$

We remark that we consider constant parameters $\nu_2, \nu_1$ for simplicity. Introducing the diagonal matrix $\bR\coloneqq\operatorname{diag}(\rho^\infty_1,\,\dots,\,\rho^\infty_n)$, this linear system may be rewritten in compact form as
\begin{equation}
	\frac{d\tilde{\bmm}}{dt}=\left[\chi\left(\bR^{-1}\bP\bR-\bI\right)+\sigma\nu_2\bR-\gamma\nu_1\bI\right]\tilde{\bmm},
	\label{eq:tilde_m}
\end{equation}
where $\tilde{\bmm}\coloneqq(\tilde{m}_1,\,\dots,\,\tilde{m}_n)$, whence we deduce that the stability of the asymptotic state $\bmm^\infty=\bzero$ depends on the spectral properties of the matrix
$$ \bA\coloneqq\chi\left(\bR^{-1}\bP\bR-\bI\right)+\sigma\nu_2\bR-\gamma\nu_1\bI\in\R^{n\times n}. $$

Remark that $\bA$ has the form $\bA=\bB-\bD$ with
$$ \bB\coloneqq\chi\bR^{-1}\bP\bR+\mu\nu_2\bR\in\R^{n\times n}, \qquad \bD\coloneqq\chi\bI+\gamma\nu_1\bI\in\R^{n\times n}, $$
where $\bB$, $\bD$ are both non-negative and $\bD$ is diagonal and invertible (at least for either $\chi>0$ or $\gamma\nu_1>0$). The Perron-Frobenius theory allows to state that $\tilde{\bmm}=\bzero$ is a stable equilibrium of~\eqref{eq:tilde_m} if and only if the spectral radius of the matrix $\bB\bD^{-1}$ is smaller than $1$. Conversely, if such a spectral radius is larger than $1$ then $\tilde{\bmm}=\bzero$ and $\bmm^\infty=\bzero$ are unstable. 

Let $\beta_{ij}\in\R$ be the $ij$-entry of the matrix $\bB\bD^{-1}$, then:
$$ \beta_{ij}\coloneqq\sum_{k=1}^{n}\bB_{ik}\bD^{-1}_{kj}=\frac{\chi P_{ij}\rho^\infty_j}{\chi\rho_i^\infty+\gamma\nu_1 \rho_i^\infty}+\dfrac{\sigma \nu_2 \rho_i^\infty}{\chi+\gamma \nu_1}\delta_{ij}, $$
and, from the Perron-Frobenius Theorem,
\begin{equation}
\min_{i=1,\,\dots,\,n}\sum_{j=1}^{n}\beta_{ij}\leq\cR_0\leq\max_{i=1,\,\dots,\,n}\sum_{j=1}^{n}\beta_{ij}.
	\label{eq:PF}
\end{equation}

In conclusion, we have that
\begin{equation}\label{eq:R0}
\dfrac{\chi +\sigma \nu_2 \rho_m^\infty}{\chi+\gamma \nu_1}   \le \cR_0 \le \dfrac{\chi +\sigma \nu_2 \rho_M^\infty}{\chi+\gamma \nu_1}, 
\end{equation}
being $\rho_m^\infty=\min_{i \in \cI} \rho_i^\infty$ and $\rho_M^\infty=\max_{i \in \cI} \rho_i^\infty$.

We can also determine the basic reproduction number in the controlled case: using the same estimation we find
\begin{equation}\label{eq:R0-control}
\displaystyle \min_{i \in \cI}\dfrac{\chi (1-u_i^{\chi,\infty}) +\sigma(1-u_i^{\sigma,\infty}) \nu_2 \rho_i^\infty}{\chi(1-u_i^{\chi,\infty})+\gamma \nu_1}   \le \cR_0 \le \displaystyle \max_{i \in \cI}\dfrac{\chi(1-u_i^{\chi,\infty}) +(1-u_i^{\sigma,\infty})\sigma \nu_2 \rho_i^\infty}{\chi(1-u_i^{\chi,\infty})+\gamma \nu_1}.
\end{equation}

\section{Numerical experiments}\label{sec:num}

In this section, we present several results from numerical experiments on the previously discussed models. In particular, we structure the presentation in three parts:
\begin{enumerate}
    \item Firstly, we show different aggregate trends arising from the most basic type of interactions on a graph, both in presence and absence of different control strategies. The equations we are considering here are~\eqref{eq:rho_i}, \eqref{eq:rho_im_i} and~\eqref{eq:brho_u},~\eqref{eq:momentum_u}, with control choices given by equations~\eqref{eq:u}--\eqref{k_mu}. 
    \item Next, we consider the descriptions gathering together infection and healing, again in both the controlled and the uncontrolled case. This time, the equations we focus on are~\eqref{eq:rho_i},~\eqref{eq:rho_im_i_new} and~\eqref{eq:brho_u},~\eqref{eq:momentum_u_new}, with control choices given by equations~\eqref{eq:u_new},~\eqref{eq:u_star_new} and~\eqref{k_chi}-\eqref{k_interval}, choosing the lowest value for the penalization coefficient. Finally, we use~\eqref{eq:R0},\eqref{eq:R0-control} to estimate the basic reproduction number.
    \item We conclude the section by applying the infection-healing model to a real-world mobility scenario, employing recent census data on the northern Italy region of Lombardy.  
\end{enumerate}
All the aforementioned macroscopic equations are solved via standard numerical methods for systems of  ordinary differential equations, such as the classical fourth-order Runge-Kutta method. We always consider a fully connected graph with $N = 5$ nodes. Unless otherwise specified, we set the parameters reported in Table~\ref{tab:general} for all the numerical tests.

\begin{table}[htbp]
\centering
\fboxrule=0.4ex\relax
\fbox{\begin{minipage}{\dimexpr\textwidth-2\fboxrule-2\fboxsep}
\centering
$\displaystyle
\bP = \left(
\begin{array}{*{5}{l}}
    0.2 & 0.5 & 0.15 & 0.1 & 0.1\\
    0.2 & 0.2 & 0.45 & 0.4 & 0.2\\
    0.2 & 0.1 & 0.05 & 0.2 & 0.5\\
    0.2 & 0.1 & 0.1 & 0.15 & 0.1\\
    0.2 & 0.1 & 0.25 & 0.15 & 0.1
\end{array}\right), \quad
\brho(t = 0) = \left(
\begin{array}{l}
    0.35\\
    0.1 \\
    0.3 \\
    0.05\\
    0.2
\end{array}\right), \quad
\bmm(t = 0) = \left(
\begin{array}{l}
    2  \\
    4  \\
    0.1\\
    1\\
    1.5
\end{array}\right),
$
\\[1ex]
$N = 5, \quad q = 2, \quad \chi = 1.$
\end{minipage}}
\caption{Parameters and initial conditions common to all tests in Section~\ref{sec:num1}--\ref{sec:num2}, unless otherwise specified.}
\label{tab:general}
\end{table}

\subsection{Test 1: infection dynamics only}\label{sec:num1}

\begin{figure}[htbp]
    \hbox to \textwidth{%
    \includegraphics[width = 0.45\linewidth]{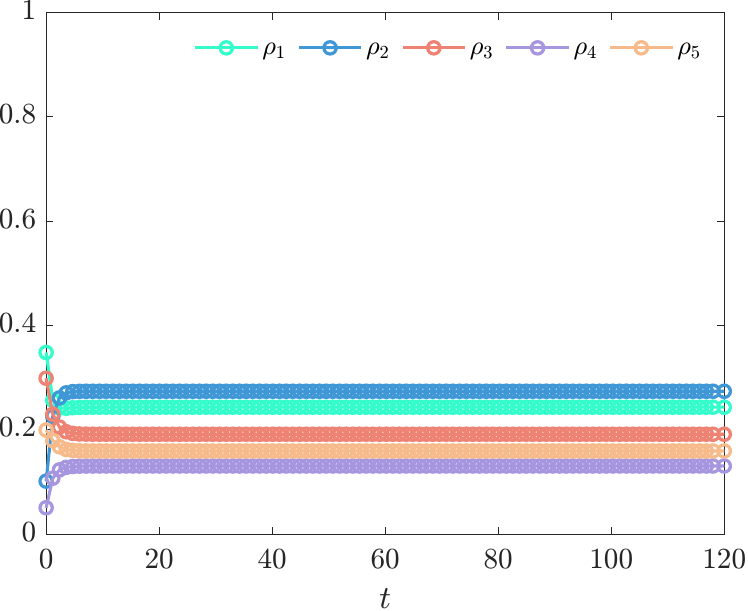}\hfil
    \includegraphics[width = 0.45\linewidth]{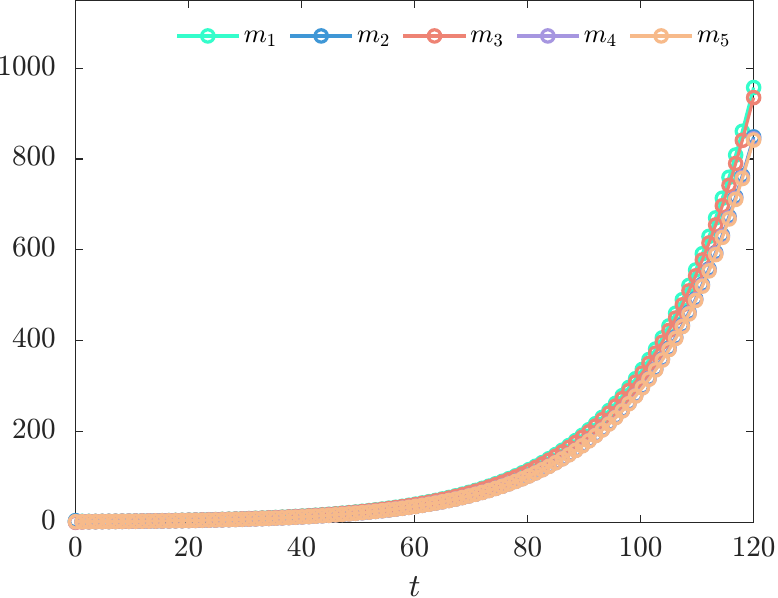}}
    \caption{Evolution in time of number of agents (left) and average viral load (right): uncontrolled case. The infection grows exponentially on the graph.}
    \label{fig:uncontrolled}
\end{figure}

\begin{figure}[htbp]
    \setbox0=\hbox{\includegraphics[width = 0.45\linewidth]{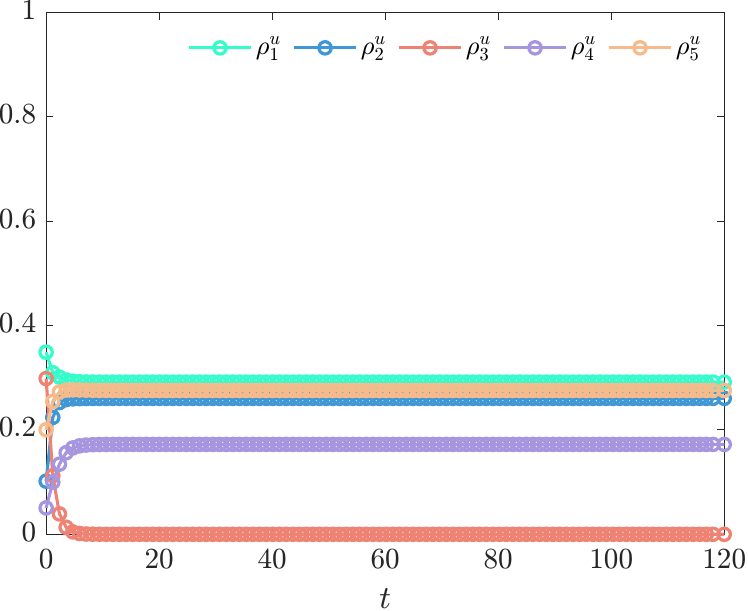}}
    \hbox to \textwidth{%
    \includegraphics[width = 0.45\linewidth]{plots/rho-mobility-only-new.pdf}\hfil
    \includegraphics[width = 0.466\linewidth, height=\ht0]{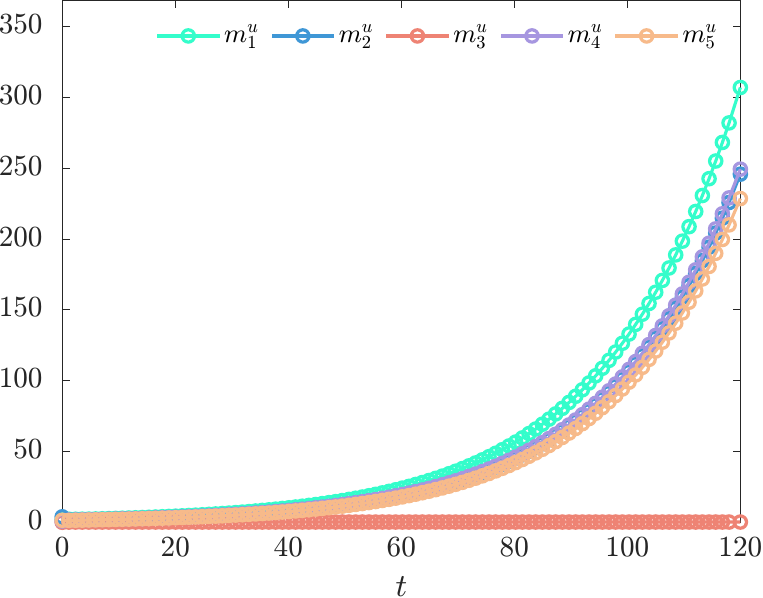}}

    \medskip
    \hbox to \textwidth{%
    \includegraphics[width = 0.45\linewidth]{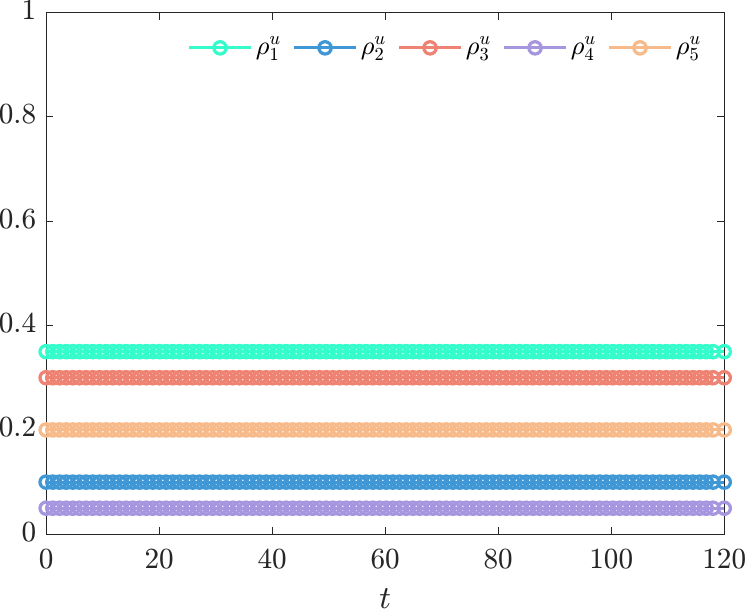}\hfil
    \includegraphics[width = 0.45\linewidth]{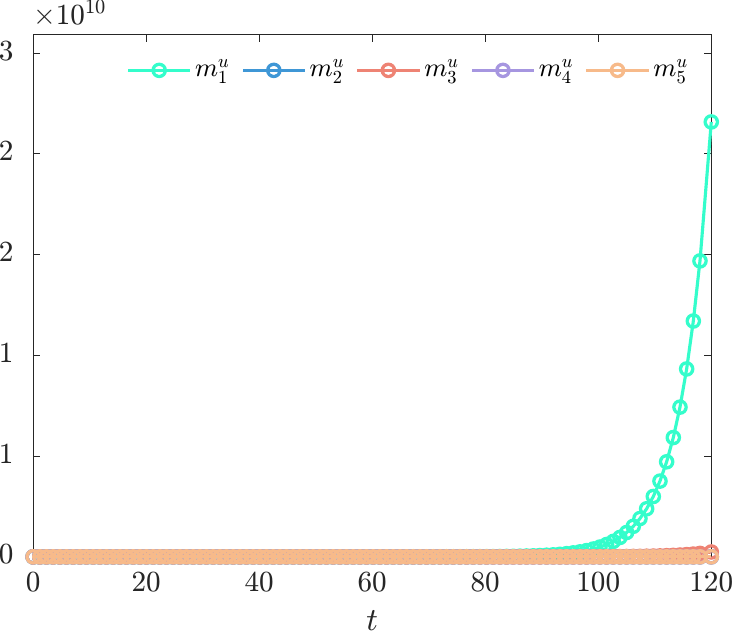}}
    \caption{Columns: evolution in time of number of agents (left) and average viral load (right). Top row: effects of controlling the agents' mobility alone. Agents leave node number 3 and distribute in the remaining vertices. The overall viral load is decreased by a factor of 3 with respect to the uncontrolled scenario. Bottom row: effects of node isolation. Without in-node interventions and preventing highly infectious agents to distribute in nodes with lower load, the infection can grow dramatically by several order of magnitudes.}
    \label{fig:mobility}
\end{figure}

\begin{figure}[htbp]
\setbox0=\hbox{\includegraphics[width = 0.45\linewidth]{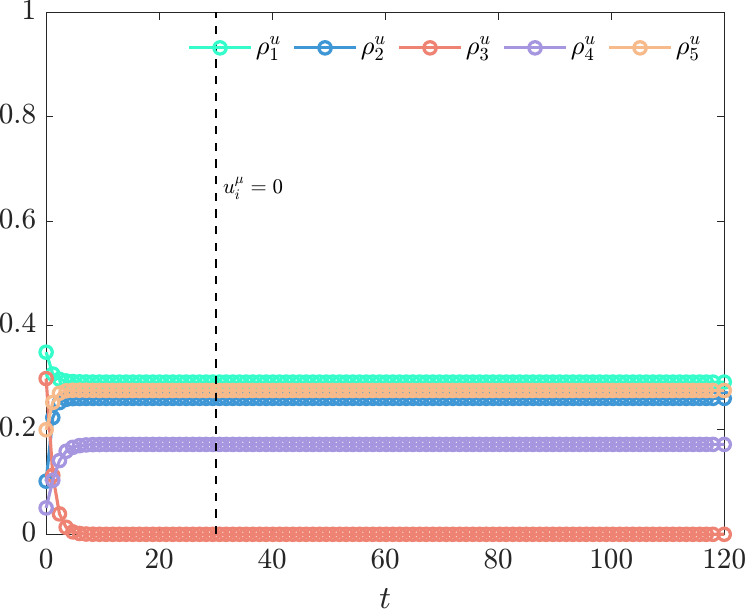}}
    \hbox to \textwidth{%
    \includegraphics[width = 0.45\linewidth]{plots/rho-mu-early-stopped-new.pdf}\hfil
    \includegraphics[width = 0.466\linewidth, height=\ht0]{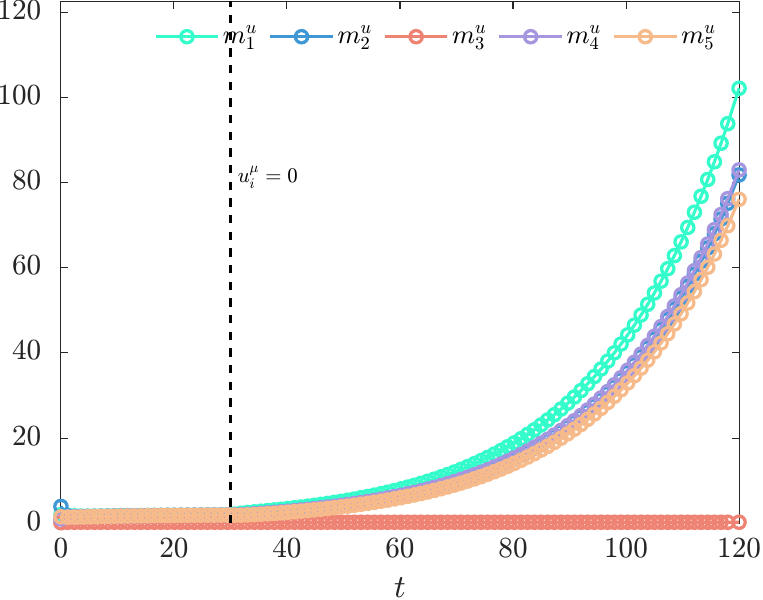}}
    
    \medskip
    \hbox to \textwidth{%
    \includegraphics[width = 0.45\linewidth]{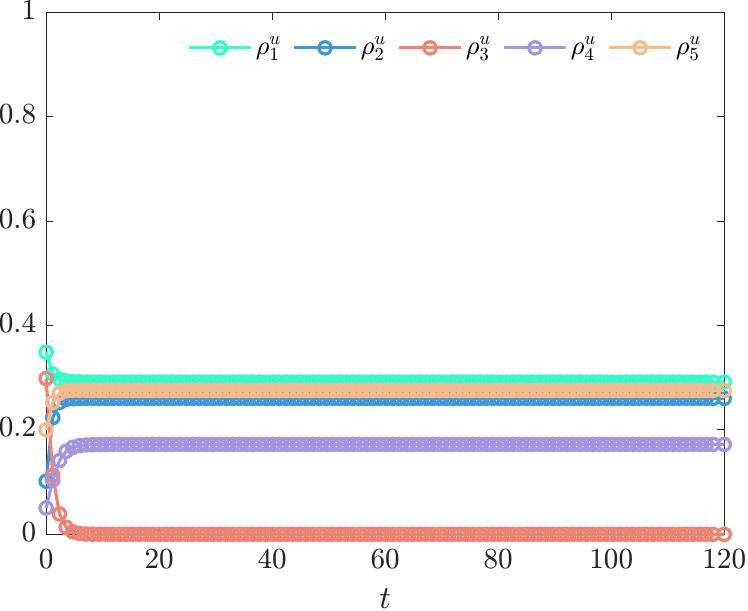}\hfil
    \includegraphics[width = 0.45\linewidth]{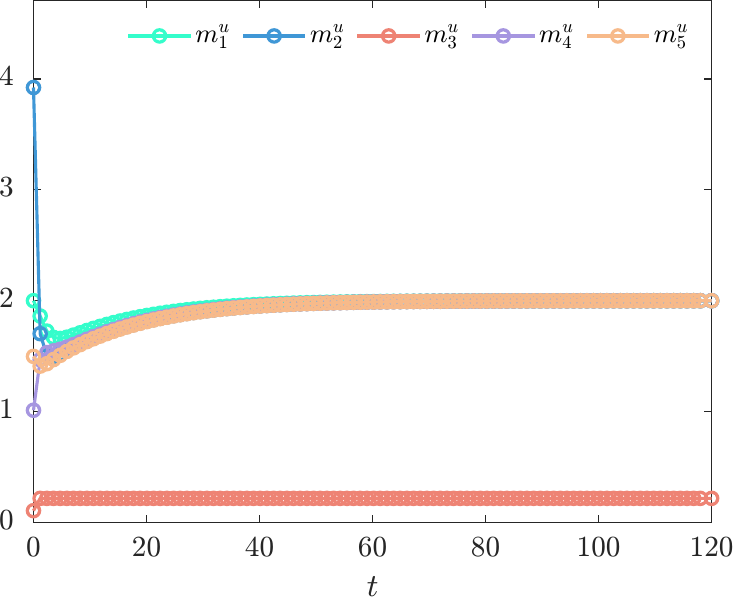}}
    \caption{Columns: evolution in time of number of agents (left) and average viral load (right). Top row: effects of partial in-node interventions. Even if limited to the early stages of infection, when it is still spreading to a comparatively low levels, control policies that slow in-node interactions between agent have a meaningful impact when compared to interventions on the mobility alone, with a reduction of about 30\% in the average viral load values. Bottom row: both in-node interactions and mobility are controlled fully. The infection ceases spreading after the average viral load reaching a value nearly three orders of magnitudes lower than the uncontrolled scenario.}
    \label{fig:interactions-and-full}
\end{figure}

We begin with the basic model described in Section~\ref{sec:basic} along with its controlled version described in Section~\ref{sec:control}. We present five different evolution scenarios, depending on the control strategy in force
\begin{enumerate}
    \item No control: this is our reference, in order to compare the effects of different intervention policies (Figure~\ref{fig:uncontrolled}).
    \item Control on mobility only: we set $u_i^\chi(t)$ according to equation~\eqref{eq:u} but we set $u_i^\mu(t) = 0$ for all $t \ge 0$ (Figure~\ref{fig:mobility}, top row).
    \item Mobility suppression: we set $u_i^\chi(t) = 1$ for all $t \ge 0$. This behavior mimics the effects of an enforced, total quarantine over the entire network. No control is set on the infection dynamics (Figure~\ref{fig:mobility}, bottom row).
    \item Short-term intervention on interactions: both dynamics are controlled until the evolution time reaches a threshold value $t = \bar t$, after which, control on interactions is suspended, i.e., $u_i^\mu (t) = 0$ for all $t > \bar t$. In our test, we set $\bar t = 30$ (Figure~\ref{fig:interactions-and-full}, top row).
    \item Full control: $u_i^\chi(t)$ and $u_i^\mu(t)$ are set according to~\eqref{eq:u} for all $t \ge 0$ (Figure~\ref{fig:interactions-and-full}, bottom row).
\end{enumerate}
For all tests, we additionally set the following parameters:
\(\label{eq:par-no-healing}
\nu_1 = (0.25,\, 0.5,\, 0.15,\, 0.2,\, 0.75)^T, \quad
\nu_2 = (0.8,\, 0.5,\, 0.75,\, 0.1,\, 0.6)^T, \quad  \mu = 1.
\)
Figure~\ref{fig:uncontrolled} shows the reference, uncontrolled case, where, at the end of the simulation period, the average viral load in each node reaches values near $10^3$. Intervening on just the mobility, rerouting infectious people, can be partially effective (see Figure~\ref{fig:mobility}, top row, exhibiting the average viral load reduced by a factor of 3); still, it is insufficient to prevent the growth of infection overall to satisfactory levels. As shown in Figure~\ref{fig:mobility}, bottom row, isolating nodes, without simultaneous in-node interventions, can not only be ineffective, but even cause more harm than good, as already remarked, e.g., in~\cite{espinoza20} and~\cite{loytosin2021}. Indeed, in-node interventions, even if for just a shorter amount of time, have a high impact on the growth rate of the average viral load: controlling the in-node interactions for just the initial 25\% of time is responsible for a nearly 33\% decrease in the infection spread, as shown in Figure~\ref{fig:interactions-and-full}, top row.

Finally, the bottom row of Figure~\ref{fig:interactions-and-full} shows the effectiveness of combining the control actions on both mobility and in-node interactions, enabling the infection spread to stop at a low average viral load value.

\subsection{Test 2: infection and healing dynamics altogether}\label{sec:num2}

Next, we focus on the model coupling the binary, in-node interaction dynamics with a healing process, described in Section~\ref{sec:healing}. In this test, we compare the uncontrolled evolution of the infection with a fully controlled scenario, where we set $u_i^\chi$ and $u_i^\sigma$ as prescribed by equation~\eqref{eq:u_star_new}. In particular, in addition to the parameters and data in Table~\ref{tab:general}, we also set
\(\label{eq:par-healing}
\nu_1^i = 0.15, \qquad \nu_2^i = 0.9, \quad \forall\, i \in \cI, \quad
\sigma = 1, \quad \gamma = 1,
\)
the exchange parameters being chosen in order to ensure that $m_i(t) \to +\infty$ in the uncontrolled case (see Proposition~\ref{prop12}), so that we can have a sensible comparison of both dynamics.

\begin{figure}[htbp]
    \setbox0=\hbox{\includegraphics[width = 0.3\linewidth]{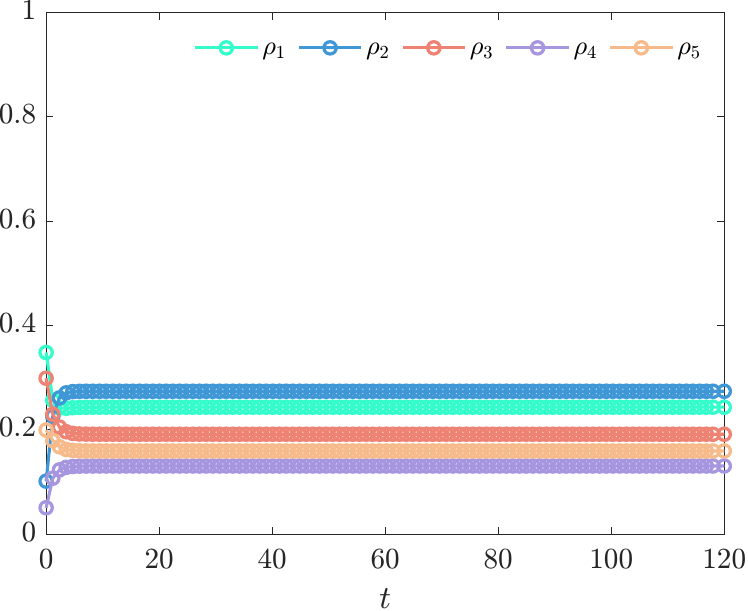}}
    \hbox to \textwidth{%
    \includegraphics[width = 0.3\linewidth]{plots/rho-uncontrolled-healing-new.pdf}\hfil
    \includegraphics[width = 0.3\linewidth, height=\ht0]{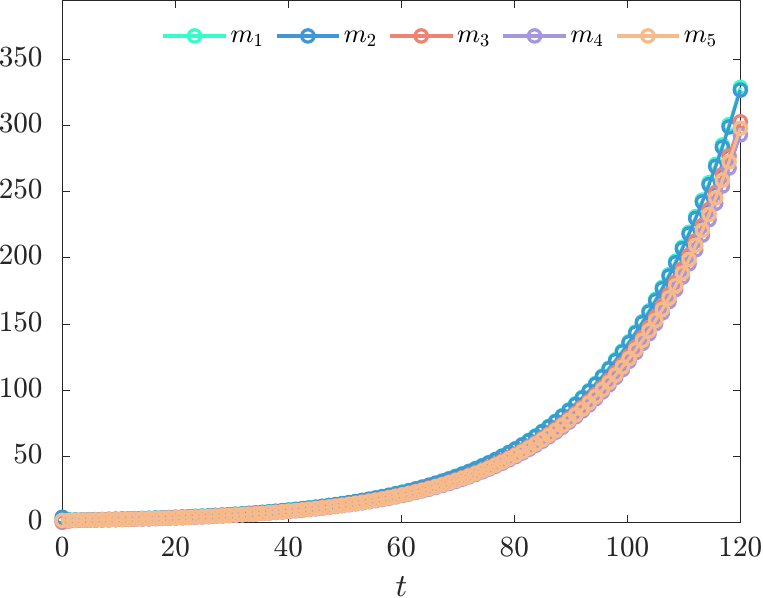}\hfil
    \includegraphics[width = 0.3\linewidth, trim=1cm 0.075cm 1cm 0.8cm, clip, height=\ht0]{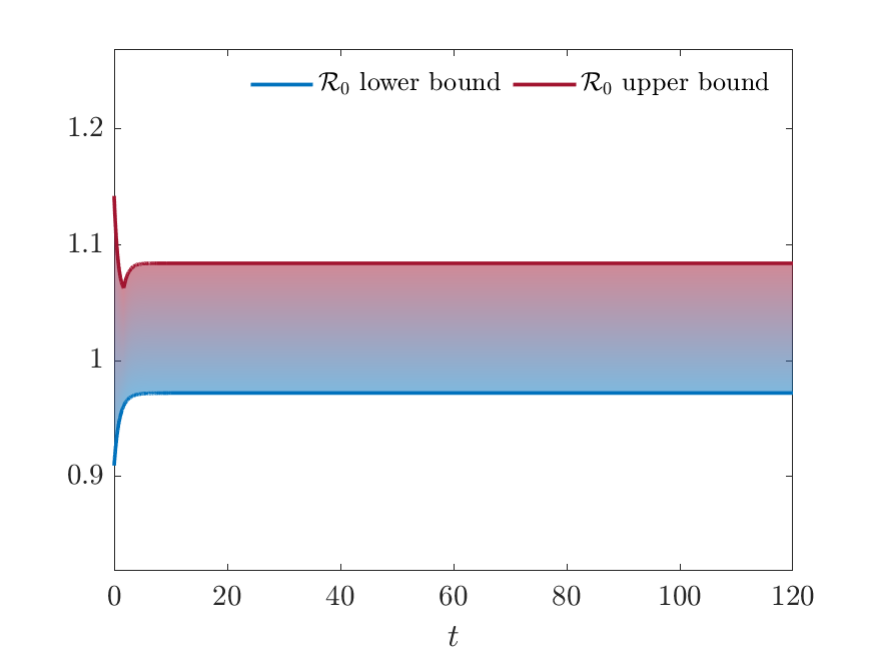}}
    
    \medskip
    \setbox0=\hbox{\includegraphics[width = 0.3\linewidth]{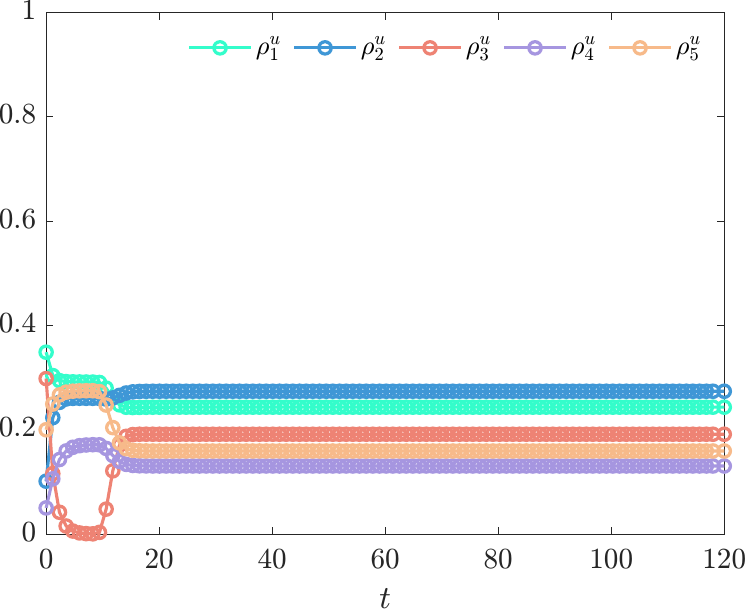}}%
    \hbox to \textwidth{%
    \includegraphics[width = 0.3\linewidth]{plots/rho-full-control-healing-new.pdf}\hfil
    \includegraphics[width = 0.3\linewidth, height=\ht0]{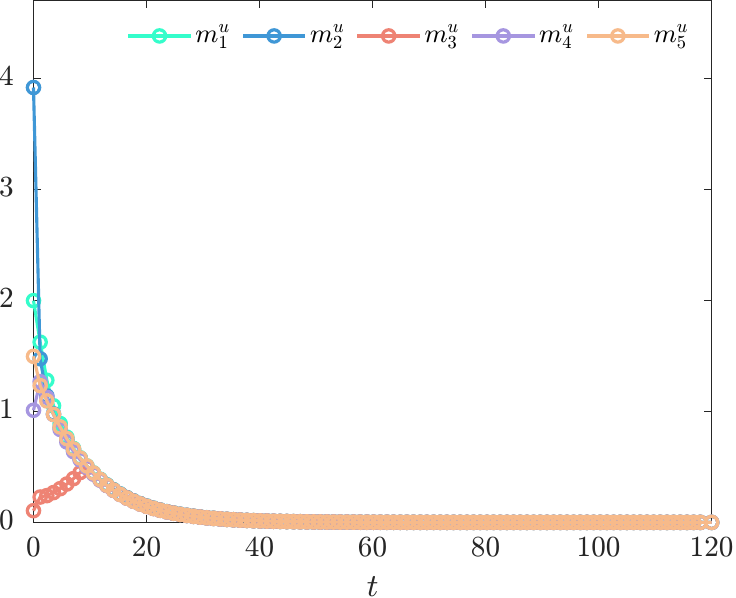}\hfil
    \includegraphics[width = 0.3\linewidth, height=\ht0, trim=1cm 0.075cm 1cm 0.8cm, clip]{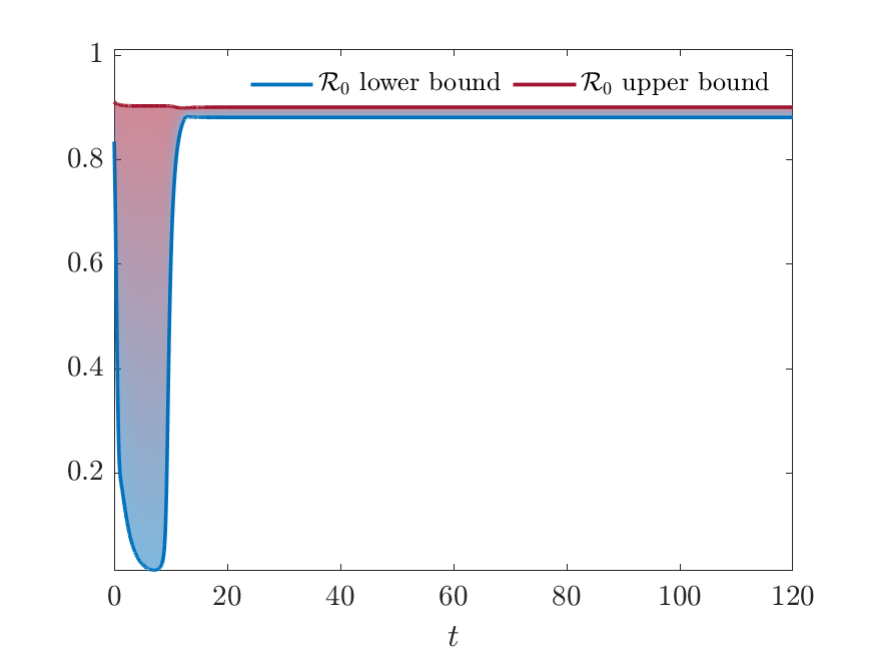}}
    \caption{Columns, from left to right: evolution in time of the number of agents, average viral load and basic reproduction number bounds. Top row: uncontrolled case, as reference. Even if slower than the more basic model without healing, we still have exponential growth of the average viral load. Bottom row: effects of full control. Infection never starts and we see the system exponentially reaching a disease-free equilibrium.}
    \label{fig:healing}
\end{figure}

In Figure~\ref{fig:healing} we report the results of the simulations: as expected, the healing process is capable alone to slow down the spread of the infection, which reaches lower values than the ones reported in its basic counterpart in Figure~\ref{fig:uncontrolled}, even with exchanging coefficients $\nu_1^i$ and $\nu_2^i$ more prone to faster dissemination. Nevertheless, the top row shows that the average viral load still grows exponentially in the uncontrolled case. This is also testified by the upper bound for the basic reproduction number being greater than 1, computed as reported in equation~\eqref{eq:R0}.

The bottom row instead shows the effects of controlling both the mobility and the in-node interactions: controlling the latter has so much relative importance, in this example, that the evolution of the number of agents is barely affected and only in the initial stages of the simulation. On the other hand, the control strategy is highly effective on the contagion dynamics, as the average viral load vanishes rapidly, as testified again by the bounds for the basic reproduction number, both below the critical value 1. These results are in agreement with the computations presented in Section~\ref{sec:healing}.

\begin{figure}
    \centering
    \hbox to \textwidth{%
    \hfil
    \includegraphics[width = 0.45\linewidth]{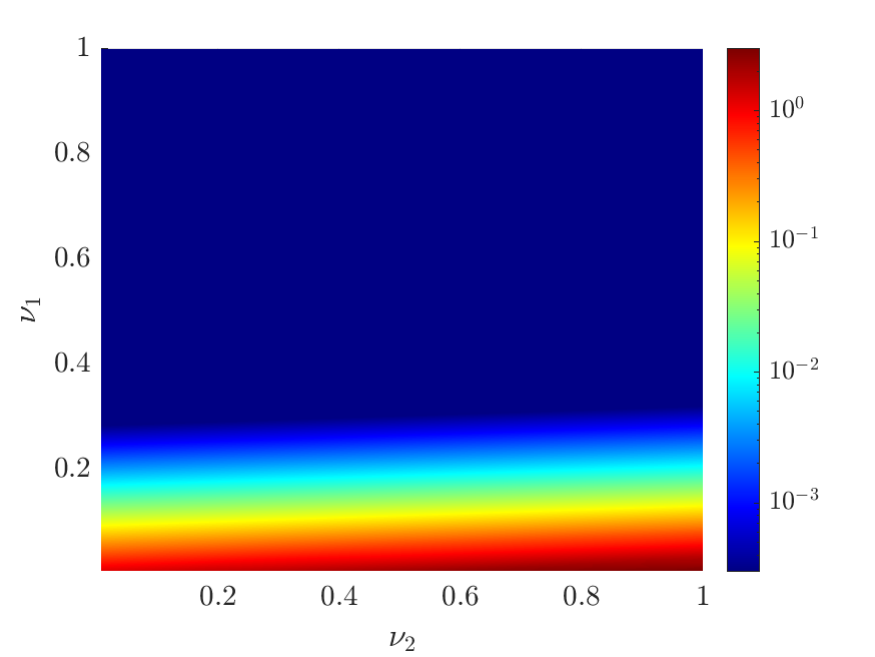}\hfil
    \includegraphics[width = 0.45\linewidth]{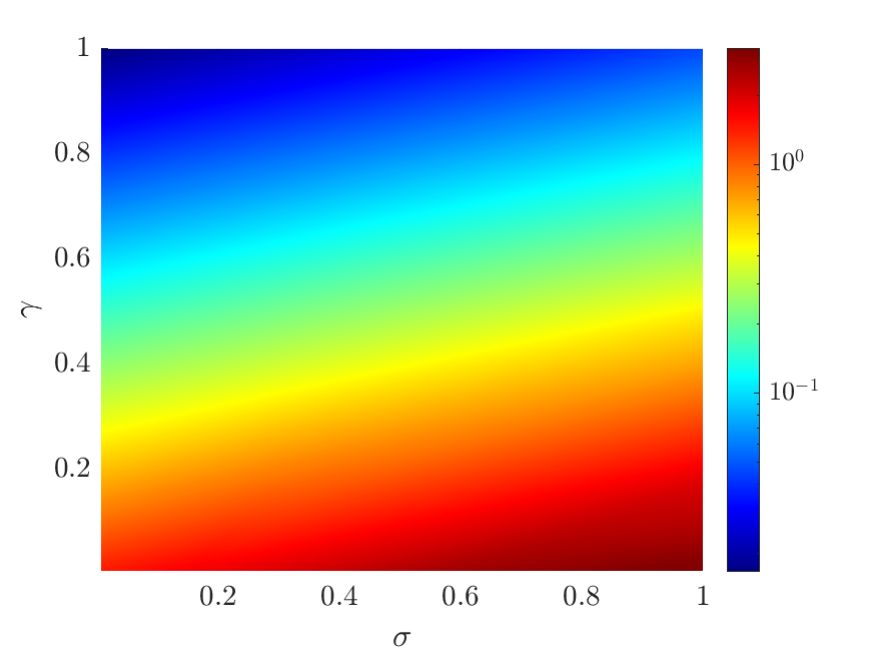}\hfil}
    \caption{Effects of parameters choice on control's efficacy. Left: largest value of $m_i$ at time $T = 30$ as a function of $\nu_1^i$ and $\nu_2^i$. Right: largest value of $m_i$ at time $T = 30$ as a function of $\gamma$ and $\sigma$. When $\nu_1^i > \nu_2^i$, we have that $m_i \to 0$, and the same when $\gamma\nu_1^i > \sigma\nu_2^i$, in accordance with Proposition \ref{prop11}. The control is less effective only when there is a relatively high ratio $(\sigma\nu_2)/(\gamma\nu_1)$.}
    \label{fig:sensitivity}
\end{figure}

In Figure \ref{fig:sensitivity} we show how different choices of parameters can affect the control. We plot the highest value of the mean viral load at time $T = 30$ as a function of pairs of $\nu_1^i$, $\nu_2^i$  and of $\gamma$ and $\sigma$. The free parameters are kept consistent with the choice~\eqref{eq:par-healing}. We see that the combinations that result in a higher value of $\rho_i^c$ are more effective in lowering the mean viral load in the nodes, in accordance with Proposition~\ref{prop11}.

\subsection{Application to a real-world mobility scenario}\label{sec:lombardy}

We conclude the presentation of numerical experiments with an application to a real-world scenario. We chose to focus on data on mobility only for two reasons: first of all, even if we are presenting here models on agents exchanging viral load, the modeling framework is flexible enough to be oriented to different (or more general) binary interactions between agents on an underlying network, indeed being the graph and its associated transition matrix the critical components of our framework. Moreover, data on viral loads often present some degree of criticality\footnote{See for example~\cite{kelley07, petrie13, hay21, killingo17} and references therein.} (under-representation of low values of viral loads associated with little to no symptoms, lack of extensive testings in non-hospitalized patients, variable fitness of the carrying pathogen over time, \ldots), making it extremely difficult to calibrate a multi-agent mathematical model on them (the interested reader may refer to e.g.,~\cite{bondesan24} for a very recent data-driven approach with an underlying kinetic framework). 

\begin{figure}
    \centering
    \includegraphics[width=0.33\linewidth]{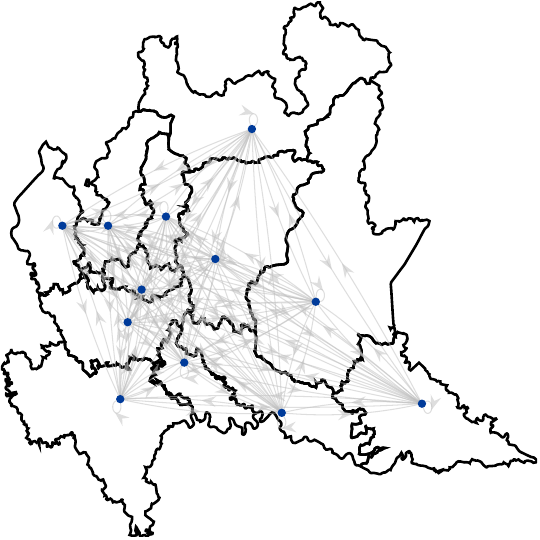}
    \caption{Fully connected, unweighted mobility network among all the provinces in Lombardy, Italy.}
    \label{fig:lombardy-basic}
\end{figure}

In order to keep the dataset relevant to our modeling setting (that is, a relatively large network, highly populated but still small enough to allow daily external mobility, so larger than a metropolitan area but smaller than a country) and also relevant for the recent SARS-CoV-2 pandemic, we chose data about the mobility habits of inhabitants of the region Lombardy, in northern Italy in 2016 (pre-pandemic)\footnote{Mobility data available at \url{https://www.dati.lombardia.it/browse?tags=80MatriceP}, last visited: 2024/08/29}. The spatial representation of the network is presented for reference in Figure~\ref{fig:lombardy-basic}. For what concerns the initial number of agents in each node, we considered the official country's census data\footnote{Population data available at\url{http://dati-censimentopopolazione.istat.it/Index.aspx}, last visited: 2024/08/29. Since we are referring to 2016 for mobility data and census is conducted once every ten years, we averaged data of 2011 and 2021. The numerical experiment results did not change significantly considering either 2011 or 2021 data, since they did not differ significantly.}. The corresponding transition matrix and initial data and parameters are reported in Table~\ref{tab:lombardy}.

\begin{table}[hbtp]
\scriptsize
\centering
\fboxrule=0.4ex\relax
\fbox{\begin{minipage}{\dimexpr\textwidth-2\fboxrule-2\fboxsep}\centering
$\displaystyle
\bP = \left(
\begin{tblr}{colspec={S[table-format=1.4]S[table-format=1.4]S[table-format=1.4]S[table-format=1.4]S[table-format=1.4]S[table-format=1.4]S[table-format=1.4]S[table-format=1.4]S[table-format=1.4]S[table-format=1.4]S[table-format=1.4]S[table-format=1.4]}, colsep=1ex}
 0.8757873534506977 & 0.03526595143954377 & 0.0054794387934947435 & 0.02158985767625072 & 0.03736936844586586 & 0.002302956167961513 & 0.011334416300349221 & 0.013701264437385683 & 0.0029183476450934482 & 0.002845914850987816 & 0.0036488799229510774 & 0.0032052497810947764\\
0.04635026989878317 & 0.9257936902391382 & 0.0014357475896173866 & 0.03174345308918152 & 0.002270751201990444 & 0.003594862553021278 & 0.003661559967382806 & 0.0050014536139418265 & 0.034928299968198774 & 0.002278835850734442 & 0.0031427242639248897 & 0.002446923335611406\\
0.002295776667751961 & 0.0005129948376395526 & 0.8154610223024258 & 0.0006910638094203594 & 0.034635103119906446 & 0.001548916939116723 & 0.03263588465510188 & 0.007767637575529528 & 0.0006900522857913599 & 0.0017705852685289778 & 0.015366268112956225 & 0.024656214415861424\\
0.007162945108250052 & 0.008473783507703942 & 0.0005042923032038033 & 0.8529198120471422 & 0.00044919089214553526 & 0.041738089168026835 & 0.0008482608670884448 & 0.004101798820922688 & 0.01812697046679227 & 0.002174689338494364 & 0.00032355456510421775 & 0.0004433451163952441\\
0.010160987217152676 & 0.000519747424712148 & 0.020542885462587812 & 0.00042586302410433605 & 0.7816261016009323 & 0.0005725317202001148 & 0.03039899767048894 & 0.0038960247334306314 & 0.00029779176703147295 & 0.0006259074183011227 & 0.01616512700820268 & 0.0012645042446099597\\
0.00044570353234450126 & 0.000573603881857499 & 0.0005527865149536564 & 0.023921110835510392 & 0.0004542106282002642 & 0.7456751566716444 & 0.0004138822564682702 & 0.010470233241153605 & 0.00043709487768829396 & 0.012195850902711174 & 0.0002701413218153193 & 0.0007310066719944785\\
0.008474134237854567 & 0.0020595765975844737 & 0.05198589677940409 & 0.0018829004382285772 & 0.08172277492366109 & 0.0012378704470828133 & 0.7104608385707953 & 0.044379804233793264 & 0.0013092858347657729 & 0.0032717205435243 & 0.002271165866704664 & 0.008437447792351163\\
0.0446185868510531 & 0.013123664706978934 & 0.053052430159367024 & 0.04106804676763683 & 0.047090653451488364 & 0.17111049006913687 & 0.1978074094567671 & 0.8672917369740902 & 0.0035986247859183156 & 0.10563301428593681 & 0.004269616682098683 & 0.10645402105185095\\
0.0005365037466188645 & 0.010815293520094473 & 0.000351684818746578 & 0.02149587688610248 & 0.0001434776436996389 & 0.0008074804316553985 & 0.0004333429105568558 & 0.0003645152123566003 & 0.9349275894464925 & 0.0009036078108297096 & 0.00023193648008865394 & 0.0005333297950915316\\
0.0012880928296168286 & 0.0008864002328547886 & 0.0016211972361365747 & 0.0028965791176721973 & 0.0009985689223157095 & 0.027892767776822522 & 0.002062754282659092 & 0.015529675735539672 & 0.0011191519702831402 & 0.8647788918483131 & 0.0005719695989734889 & 0.001587329172944043\\
0.0005778134784514185 & 0.000459143073049482 & 0.005893878831690055 & 0.00018553162048871004 & 0.009780314082241892 & 0.000202691408857945 & 0.0005387291445374318 & 0.00018473638480129612 & 0.0001500703598885761 & 0.0001647479386678979 & 0.9524046220972929 & 0.0002855433643589506\\
0.0023018329814251818 & 0.0015161505388427983 & 0.04311873920837233 & 0.0011799046882619878 & 0.003459485087552358 & 0.003316186646473646 & 0.009403923917804436 & 0.02731111903705499 & 0.0014967205920561203 & 0.003356233942970127 & 0.0013339940798873212 & 0.849955085257836\\
\end{tblr}\right)$, \enspace
$\displaystyle
\brho(t = 0) = \left(
\begin{tblr}{colspec={S[table-format=1.4]}, colsep=1ex}
0.111429568301365\\
0.126797034990562\\
0.0361007484289710\\
0.0601448912068948\\
0.0340388722947419\\
0.0229591510831975\\
0.0870627833902669\\
0.318267453990158\\
0.0413704681415706\\
0.0544775057762816\\
0.0183028026195141\\
0.0890487197764766
\end{tblr}
\right)$
\\[1ex]
$\displaystyle
N = 12, \quad
m_i(t=0) =
\begin{cases}
    1/6 & \text{if } i \neq 6\\
    6 &  \text{otherwise}
\end{cases}, \qquad
\nu_1^i = 0.15, \qquad \nu_2^i = 0.9, \quad \forall\, i \in \cI$
\end{minipage}}
\caption{Parameters associated to the Lombardy mobility network and control problem, along with parameters associated to viral load.}
\label{tab:lombardy}
\end{table}

\begin{figure}[htbp]
    \setbox260=\hbox{\includegraphics[width=0.45\textwidth]{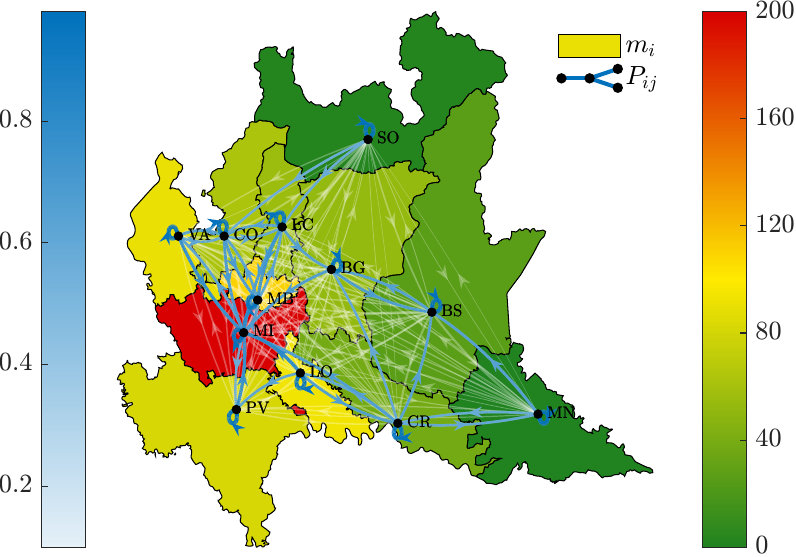}}
    \setbox270=\hbox{\includegraphics[width=0.45\textwidth]{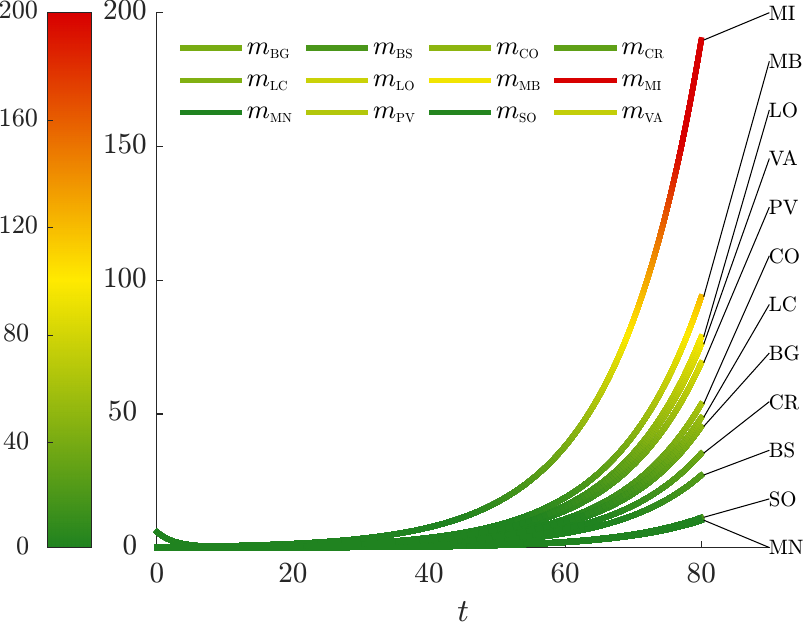}}
    \begin{minipage}[b][\ht270][t]
    {0.45\textwidth}
    \includegraphics[width = \linewidth, height=0.975\ht260]{plots/lombardy-uncontrolled.pdf}
    \end{minipage}\hfill
    \begin{minipage}[b][][t]{0.45\textwidth}
    \includegraphics[width = \linewidth]{plots/lombardy-uncontrolled-means-new.pdf}
    \end{minipage}
    
    \medskip
    \setbox260=\hbox{\includegraphics[width=0.45\textwidth]{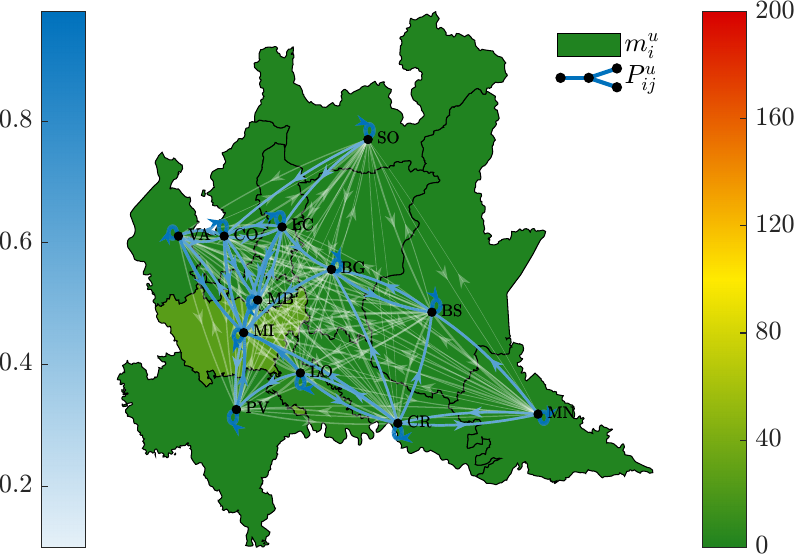}}
    \setbox270=\hbox{\includegraphics[width=0.45\textwidth]{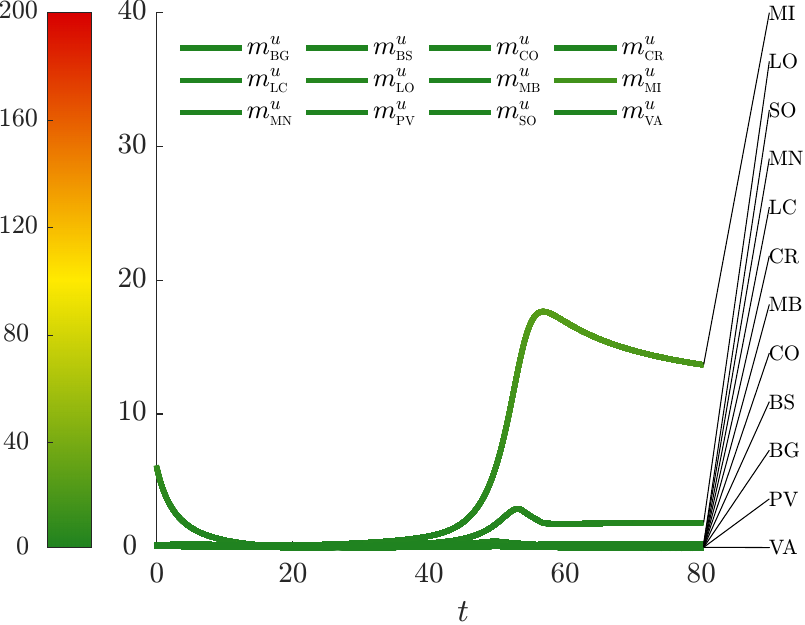}}
    \begin{minipage}[b][\ht270][t]
    {0.45\textwidth}
    \includegraphics[width = \linewidth, height=0.975\ht260]{plots/lombardy-controlled.pdf}
    \end{minipage}\hfill
    \begin{minipage}[b][][t]{0.45\textwidth}
    \includegraphics[width = \linewidth]{plots/lombardy-controlled-means-new.pdf}
    \end{minipage}
    \caption{Columns, from left to right: heat-map of the infection state along the mobility network at the end of the simulation period; evolution in time of the average viral load on the network. Top row: uncontrolled scenario, as reference. We see exponential spread of the infection, at different rates. Bottom row: effects of light in-node control. Even if the control strength was not sufficient to achieve eradication across the network, the average viral load reaches values that are more than one order of magnitude lower than the reference case, thus obtaining a satisfactory result.}
    \label{fig:lombardy}
\end{figure}

In Figure~\ref{fig:lombardy} we report the results of the numerical experiment. We consider again the infection-healing modeling framework of Section~\ref{sec:healing}, as we did in Section~\ref{sec:num2}, but this time we consider much less restrictive in-node control actions. Indeed, we know from Proposition~\ref{prop12} that we can achieve eradication over the whole network if we set $u_i^\sigma(t) \ge 1 - \rho_i^c$ for all $t \ge 0$. This, however, might imply a very high associated cost, especially in those cases when some nonzero viral load is present along the network but its value is still negligible.

In the bottom row of Figure~\ref{fig:lombardy} we show that even if we do not achieve complete eradication, we still are able to lower the average viral load by more than one order of magnitude by relaxing the lower bound on $u_i^\sigma$ to be
\[
u_i^\sigma (t) = \min(\max(0, \tilde k_i^\sigma\times 2.5\cdot 10^{-5}), 1), \qquad \tilde k_i^\sigma = (\rho_i(0)m_i(0))^q \sigma \nu_2^i,
\]
with a substantial decrease in the associated penalization coefficient (and therefore its associated cost).

\section{Conclusions}

In the present work, we have proposed the optimal control of a kinetic model describing social interactions on a graph. The basic kinetic model, that was proposed in~\cite{loytosin2021}, describes agents migrating on the nodes of a graph and exchanging a physical quantity $v>0$ as a consequence of binary interactions. Here, the control, which is exerted on the microscopic mechanisms aims at minimizing the macroscopic average of the aforementioned physical quantity $v$. In~\cite{loytosin2021}, the kinetic model describes the spread of an infectious disease on a graph: as such, the positive quantity $v$ represents the viral-load of the potentially infected individual. The binary exchange rules are linear ones, and the exchange and migration mechanisms are stochastic independent. As a consequence, in the present work we have implemented two different controls on the two mechanisms in order to minimize the quantity related to the macroscopic average viral load. In order to do this we have applied the method used in~\cite{albi2021}. Specifically, we have chosen to minimize the average $\rho_i m_i$ in each node weighted by the mass. In fact, we have shown that controlling $\rho_i m_i$ is effective as either controlling the mean $m_i$ in each node or controlling the global average $m$, but less expensive. This is due to the fact that implementing the same control everywhere (control on $m$), or controlling the same amount of mean viral-load $m_i$ regardless of the population quantity $\rho_i$, corresponds to exerting an excessive control for obtaining the same result. The best success that can be reached by this controlled model is to stop the increase of the infection, but eradication cannot be obtained, unless there are \emph{a priori} natural conditions such as $\nu_2<\nu_1$. This is due to the interaction process, that includes the infection and healing within the same process.

As a consequence, we have proposed an adaptation of the model introduced in~\cite{loytosin2021}. As the key point in the model proposed in~\cite{loytosin2021} is that it prescribes simultaneous infection and healing within the same binary interaction rule, we have split the interaction into an infection process (due to a binary interaction) and a healing one (due to an autonomous process). Then, the control has been implemented only on the migration and the infection processes. This has allowed to show that a sufficient control allows to reach eradication, even if the \emph{a priori} conditions would not allow it.  

Overall, the proposed controlled models allow to test the effect of each control strategy alone as well as the interplay of the simultaneous controlling strategies. As a consequence, similar results to the ones of previous studies~\cite{espinoza20} have been obtained. Moreover, it has been possible to highlight the drawbacks of the present control strategy, that prescribes a control that may be too high and persisting in time. This happens because of the choice~\eqref{def:psi} that aims at obtaining complete eradication ($m_i\rightarrow 0$ in each node). This reminds, for example, of the quarantine methods through PCR tests adopted during the COVID-19 pandemics: as the sensitivity of those tests was too high, then also recovered and not infectious individuals were isolated, because some viral load was still measured by the swab. Even though the PCR test and the present control are different, because the first one is a microscopic (individual) control, while the second one is a macroscopic (population level) control, this suggests, as a future perspective, that the control may be improved by adapting the cost function $\psi$ in order to demand that the average viral load is below a certain, even strictly positive, threshold, instead of requiring complete eradication. This may pose some challenges as $\psi'$ could change sign instead of being always positive.
Moreover, as a possible future work, we aim at integrating this viral-load modeling approach on a graph with a compartmental one, as done for example in~\cite{RDMNLMM,RDMNLAT2}, with the addition of control strategies.

We remark that we have presented the `less realistic' model of Section~\ref{sec:basic}, because it has allowed to introduce in~\cite{loytosin2021} the concept of exchanging viral-load allowing to characterize the state of the individual with respect to the disease without considering the epidemic compartments. In the present work, this starting model has allowed to face some difficulties in the definition of the control problem. Moreover, the simple linear interaction rules~\eqref{eq:update-v} of the model of Section~\ref{sec:basic} have proved to be effective in reproducing real phenomena of wealth exchange~\cite{cordier2005JSP} and opinion exchange~\cite{toscani2006CMS}. As a consequence, our framework could be adapted in order to consider the control in other phenomena of interest of social exchange on a graph or in presence of migrating subpopulations~\cite{bisi2022PhTB,bisi2024Physd}.

\paragraph*{Author contributions}
J.F. and N.L. conceived the project, wrote the paper and made figures.

\paragraph*{Credit authorship contribution statement}
Nadia Loy and Jonathan Franceschi: Investigation, Conceptualization, Writing -- original draft, Visualization, Writing -- review and editing.

\paragraph*{Acknowledgments}
This work has been written within the activities of the GNFM group of INdAM. N.L. and J.F. acknowledge the support of GNFM through the Research Project Progetto GNFM 2023:
``Controllo ottimo di equazioni cinetiche su grafo'' CUP E53C22001930001. J.F. and N.L. acknowledge the support of the Italian Ministry of University and Research (MUR) through the PRIN-2020  project (No. 2020JLWP23) ``Integrated Mathematical Approaches to Socio-Epidemiological Dynamics''. N.L. acknowledges the support of the Italian Ministry of University and Research (MUR) through
the grant PRIN2022-PNRR project (P2022Z7ZAJ) ``A Unitary Mathematical Framework for Modelling Muscular
Dystrophies'' (CUP: E53D23018070001) funded by the European Union - Next Generation EU.

\paragraph*{Declaration of Competing Interest}
The authors declare that they have no known competing financial interests or personal relationships that could have appeared to influence the work reported in this paper.

\paragraph*{Data availability}
Data will be made available on request.

\appendix

\section{Proof of Proposition~\ref{prop:metzler}}\label{appendix:proof}

Let $\tilde t \coloneqq \max_i \{ t_1^i, t_2^i \}$. From assumption~\ref{cond:mi} we immediately deduce the eventual non-decreasing behavior of the average. Indeed, if we rewrite~\eqref{eq:mean_u} in integral form we have
\[
\begin{aligned}
m_i(t) &= m_i(\tilde t) + \chi \int_{\tilde t}^t \biggl(\sum_{\substack{i\in \cI}} \Puij(s) \frac{\rho_j(s)}{\rho_i(s)} ( m_j(s)  - m_i(s)) \biggr)\, ds + \mu\int_{\tilde t}^t (1 - u_i^\mu(s)) (\nu_2^i - \nu_1^i) \rho_i(s) m_i(s) \, ds\\
&\ge \mu\int_{\tilde t}^t (1 - u_i^\mu(s)) (\nu_2^i - \nu_1^i) \rho_i(s) m_i(s) \, ds,
\end{aligned}
\]
from which we deduce that
\(\label{eq:media-monotona}
\dt m_i \ge \mu(1 - u_i^\mu) (\nu_2^i - \nu_1^i) \rho_i m_i \ge 0,
\)
for all $t \ge \tilde t$. An analogous computation leveraging the first inequality in~\eqref{eq:bazooka-rhoimi} of assumption~\ref{cond:rhoimi} implies also that each $\rho_i m_i$ is non-decreasing for all $t \ge \tilde t$ and all $i \in \cI$.

Next, we remark that
\(\label{eq:sum_increasing}
\dt \sum_{i\in\cI} \rho_i m_i = \sum_{i\in\cI} \mu(1 - u_i^\mu) (\nu_2^i - \nu_1^i) \rho_i^2 m_i \ge 0
\)
for all $t \ge 0$ since $\nu_2^i \ge \nu_1^i$ for all $i \in \cI$, so that the sum is monotonically non-decreasing.

We proceed by proving the boundedness of $\sum_i \rho_i m_i$ for all $t \ge \tilde t$. Let us suppose that
\[
\sum_{i\in\cI} \rho_i m_i < \sum_{i\in\cI} m_i(0)
\]
for all $t \ge \tilde t$. Then, thanks to the choice~\eqref{eq:bar_u_conv}, we have that $1-u_i^\mu >0$, and both points of the claim follow straightforwardly in view of the monotonic behavior. Otherwise, we remark that equation~\eqref{eq:sum_increasing} implies that there exists $\bar t \ge \tilde t$ such that
\[
\sum_{i\in\cI} \rho_i m_i \ge \sum_{i\in\cI} m_i(0)
\]
for all $t \ge \bar t$. We define the following, time-dependent, sets of indices:
\(
\cI^+ \coloneqq \{ i \in \cI \mid m_i(t) \ge m_i(0)\}, \qquad \cI^- \coloneqq \cI \setminus \cI^+.
\)

We remark that in view of inequality~\eqref{eq:media-monotona}, we have either one of the following conditions
\begin{enumerate}
    \item $i \in \cI^-$ for all $t \ge \bar t$;
    \item There exists $t_i^\star > \bar t$ such that
    \[
    \begin{cases}
        i \in \cI^- & \text{ if } \bar t \le t < t_i^\star\\
        i \in \cI^+ & \text{ if } t \ge t_i^\star.
    \end{cases}
    \]
\end{enumerate}

These imply that
\[
\sum_{i \in \cI^-} \rho_i m_i(t) < \sum_{i \in \cI^-} m_i (0)
\]
for all $t \ge \bar t$, so that we are left to prove that the summation done over the set $\cI^+$ does not blow up eventually. This is granted by condition~\ref{cond:rhoimi}. Indeed, in integral form we have
\[
\begin{aligned}
\sum_{i \in \cI^+} \rho_i m_i(t) &= \sum_{i \in \cI^+} \rho_i m_i(t_i^\star) 
+ \sum_{i \in \cI^+}\chi \int_{t_i^\star}^t \biggl(\sum_{j\in \cI} \Puij(s) \rho_j(s) m_j(s)  - \rho_i(s) m_i(s) \biggr)\, ds\\
&+ \underbrace{\sum_{i \in \cI^+}\mu\int_{t_i^\star}^t (1 - u_i^\mu(s)) (\nu_2^i - \nu_1^i) \rho_i^2(s) m_i(s) \, ds}_{{} = 0}\\
&\le \sum_{i \in \cI^+} \rho_i m_i(t_i^\star) \Bigl(\frac{1}{r_i(t)} + \alpha_i \Bigr) \xrightarrow[t \to +\infty]{} \sum_{i \in \cI^+} \alpha_i \rho_i m_i(t_i^\star) < +\infty,
\end{aligned}
\]
where the under-braced integral vanishes because of the control term $u_i^\mu$ being identically $1$ for all $t \ge t_i^\star$ in view of the penalization choice~\eqref{eq:bar_u_conv} and index $i$ belonging to $\cI^+$ for all times greater than $t_i^\star$.

Since the total first moment is eventually bounded, we obtain the claim.

\printbibliography

\end{document}